\documentclass[12pt]{iopart} 

\usepackage[T1]{fontenc}
\usepackage{hyperref}
\hypersetup{colorlinks,allcolors=blue}
\usepackage{eurosym}
\usepackage{makecell}
\usepackage{amsmath}
\usepackage{amsfonts}
\usepackage{graphics}
\usepackage{graphicx}
\usepackage{diagbox}
\usepackage{multirow}
\usepackage{url}
\usepackage{verbatim}
\usepackage{bm}
\usepackage{accents}
\usepackage{multicol}
\usepackage{multirow}
\usepackage{eqparbox}
\usepackage{enumerate}
\usepackage{soul}
\usepackage{bbding}
\usepackage{tablefootnote}
\usepackage{longtable}
\usepackage{xcolor}
\usepackage{color}
\usepackage{soul}

\begin{document}

\title[]{Investigating State-of-the-Art Planning Strategies for Electric Vehicle Charging Infrastructures in Coupled Transport and Power Networks: A Comprehensive Review}

\author{Jinhao Li\textsuperscript{1}, Arlena Chew\textsuperscript{1}, and Hao Wang\textsuperscript{1,2*}}

\address{\textsuperscript{1}Department of Data Science and AI, Faculty of IT, Monash University, Melbourne, Australia \\
\textsuperscript{2}Monash Energy Institute, Monash University, Melbourne, Australia \\
\textsuperscript{*}Author to whom any correspondence should be addressed.}
\ead{hao.wang2@monash.edu}

\begin{abstract}
Electric vehicles (EVs) have emerged as a pivotal solution to reduce greenhouse gas emissions paving a pathway to net zero. As the adoption of EVs continues to grow, countries are proactively formulating systematic plans for nationwide electric vehicle charging infrastructure (EVCI) to keep pace with the accelerating shift towards EVs. This comprehensive review aims to thoroughly examine current global practices in EVCI planning and explore state-of-the-art methodologies for designing EVCI planning strategies. Despite remarkable efforts by influential players in the global EV market, such as China, the United States, and the European Union, the progress in EVCI rollout has been notably slower than anticipated in the rest of the world. This delay can be attributable to three major impediments: inadequate EVCI charging services, low utilization rates of public EVCI facilities, and the non-trivial integration of EVCI into the electric grid. These challenges are intricately linked to key stakeholders in the EVCI planning problem within the context of coupled traffic and grid networks. These stakeholders include EV drivers, transport system operators, and electric grid operators. In addition, various applicable charging technologies further complicate this planning task. This review dissects the interests of these stakeholders, clarifying their respective roles and expectations in the context of EVCI planning. This review also provides insights into level 1, 2, and 3 chargers with explorations of their applications in different geographical locations for diverse EV charging patterns. Finally, a thorough review of node-based and flow-based approaches to EV planning is presented. The modeling of placing charging stations is broadly categorized into set coverage, maximum coverage, flow-capturing, and flow-refueling location models. In conclusion, this review identifies several research gaps, including the dynamic modeling of EV charging demand and the coordination of vehicle electrification with grid decarbonization. This paper calls for further contributions to bridge these gaps and drive the advancement of EVCI planning. 
\end{abstract}


\section*{Acronym}
\begin{tabular}{ll}
    GHG & Greenhouse gas \\
    EV & Electric vehicle \\
    ICE & Internal combustion engine \\ 
    IEA & International Energy Agency \\ 
    EU & European Union \\ 
    EVCI & Electric vehicle charging infrastructure \\ 
    DC & Direct current \\ 
    AC & Alternative current \\ 
    O-D & Origin-destination \\ 
    SCLM & Set covering location model \\ 
    MCLM & Maximum covering location model \\ 
    FCLM & Flow-capturing location model \\
    FRLM & Flow-refueling location model \\ 
    CFRLM & Capacitated flow-refueling location model \\ 
    BPR & Bureau of Public Roads \\ 
    GA & Genetic algorithm \\ 
    PSO & Particle swarm optimization \\ 
    B\&B & Branch-and-bound \\ 
    TN & Traffic Network \\ 
    PDN & Power distribution network \\ 
\end{tabular}

\section{Introduction} \label{sec:intro}
To reach the goal outlined by the Paris Agreement of mitigating global warming well below $2$ degrees Celsius~\cite{meinshausen2022_ParisAgreement}, nations have devised their policies to curb greenhouse gas (GHG) emissions based on the nationally determined contributions~\cite{fankhauser2022_netZero}, e.g., Germany in 2045~\cite{de2021_netZero}, the U.S. in 2050~\cite{us2021_netZero}, China in 2060~\cite{china2021_netZero}, and India in 2070~\cite{india2022_netZero} for the net-zero targets~\cite{UN2023_netZero}. Among various energy-related sectors, the transport sector still heavily relies on fossil fuel resources~\cite{wu2018}, contributing significantly to GHG emissions, which account for about $15\%$ in 2021~\cite{metais2022review} and is estimated to rise to $77\%$ by 2055~\cite{gupta2021review}. Therefore, urgent and effective efforts are imperative to decarbonize the transport sector. Electric vehicles (EVs) have emerged as the most promising and practical solution to replace internal combustion engine (ICE) vehicles, since their low-carbon mobility, also known as electromobility, aligns with the net-zero transition~\cite{lamonaca2022review}. Moreover, benefiting from the emerging battery technologies, the levelized cost of EVs is expected to be lower than that of ICE vehicles by 2025~\cite{soulopoulos2017}, making EVs a more cost-effective and viable choice for electrifying and decarbonizing the transport sector. The International Energy Agency (IEA) has reported an exponential growth of EVs on road transport, with sales exceeding $10$ million in 2022~\cite{IEA2023_EVOutlook}, and an accumulated on-road EV count of $125$ million by 2023~\cite{chen2020review_UK}. Furthermore, the global market share of EVs is projected to proliferate to $35\%$ in 2030, potentially leading to the avoidance of approximately $700$ million metric tons of carbon emissions~\cite{IEA2023_EVOutlook}.

While increasing EV adoption is poised to significantly reduce carbon emissions within the transport sector, EV drivers have expressed concerns regarding the availability of charging stations due to the well-known problem of range anxiety and the inconvenience of finding vacant chargers~\cite{ji2018review_China}. It is noteworthy that charging convenience plays a pivotal role in affecting individuals' decisions to purchase an EV~\cite{hoed2013_Netherlands,wolbertus2016_Netherlands}. Therefore, to maintain the momentum of EV uptake for sustainable transportation electrification, providing ample EV charging infrastructure (EVCI) is indispensable to cater to the escalating charging needs of EV drivers. Furthermore, research by the International Council on Clean Transportation has unveiled significant correlation between the deployment of EV charging stations and the market share of EVs~\cite{hall2017}. Motivated by the strong need to propel the supply of charging infrastructure, city planners have been funded by governments to promote the expansion of EVCI. For example, France has announced the Energy Transition for Green Growth Act, aiming to build $7$ million public and private chargers by 2030, representing a minimum cost of approximately \euro$2$ billion~\cite{fr2019_energyAct}. In the United States, a one-billion-dollar investment plan is underway to bridge the gap in public EVCI by 2020~\cite{nicholas2019_US}. The Australian Government has doubled its funding for its ``Driving The Nation'' project to AU\textdollar$500$ million to build the backbone of a national EV charging network~\cite{au2023_DriveTheNation}.

Despite the above governmental policies incentivizing charging station deployment, one major concern revolves around inefficient EVCI planning. Inefficiency may compromise the interests of several stakeholders, including infrastructure developers (such as city planners) and infrastructure users (such as EV drivers), potentially slowing EV uptake. On one hand, city planners, being the primary developers in charging stations, have noticed unexpectedly low utilization rates of public EVCIs. In the Netherlands, the U.S., and China, the utilization rates are approximately $4\%$, $7\%$, and $15\%$, respectively~\cite{unterluggauer2022review,kumar2013_US,yang2020_China}, indicating that the current infrastructure planning fails to capture EV drivers' charging behaviors and does not sufficiently attract them to use public facilities. The significant underutilization of public charging resources may undermine the operator's confidence in building more EV chargers to keep up with the accelerating transition to EVs. Nevertheless, the IEA estimates that a staggering $13$ million public chargers need to be installed by 2030 to meet the projected EV charging demand of approximately $400$ TWh in 2030~\cite{IEA2022_STEPS}, indicating a substantial lack of investment~\cite{JDPower2022}. On the other hand, though EV drivers have access to vacant public EV chargers due to the high likelihood of facility availability, over half of them prefer using private chargers to charge their EVs at home~\cite{hardman2018_homeCharging}, attributed to two major reasons. First, home charging offers better accessibility~\cite{smart2015_homeCharging} as EV drivers have limited choices of public charging stations that align with their travel patterns and daily driving ranges~\cite{lee2020_homeCharging}. Second, a lack of charging resources in residential areas~\cite{schurmann2017_lackResidentialCharging}: Despite home parking consituting over $75\%$ of an EV's daily parking time~\cite{kouka2020_lackResidentialCharging}, existing public EVCI tend to be located far from residential areas, prompting EV drivers to opt for private EV chargers. Regardless of the upfront costs of installing a private charger, the prevalence of home charging emphasizes the underperformance of current EVCI.

Hence, the current inefficiency in EVCI planning has created a chicken-and-egg dilemma for charging station expansion and EV market diffusion~\cite{shi2021_chickenEggsDilemma}. The city planner hesitates to invest more in EVCIs due to the low utilization rates, while individuals are reluctant to buy an EV without easily accessible charging stations to ensure charging availability. This dilemma could significantly impede the expected EV adoption and transport decarbonization. To mitigate the challenge, effective planning for building EVCI is crucial, which should allocate available charging resources to meet EV drivers' charging demand, improve the utilization rates of charging stations, encourage EVCI developers for further investment, and provide positive signals to EV drivers to promote continuous EV adoption.

However, planning an efficient EV charging network is a complex task due to the intricate interests among stakeholders, multiple optional charging technologies for charger installation, and unclear EVCI planning strategies in large-scale traffic networks with budget considerations. To guide future directions in EVCI planning, in this paper, we aim to address the following questions:
\begin{itemize}

    \item What are the current real-world practices in building EV charging networks and what barriers are hindering the rollout of EVCIs?

    \item Who are the stakeholders involved in the specific EVCI planning problem, what are their roles and interests, and what are the advantages and disadvantages of adopting different optional charging technologies in potential facility locations?

    \item What methodological solutions have been proposed to design efficient EVCI planning strategies?

    \item What are the research gaps that require further contributions?

\end{itemize}
The first three questions have been partially explored in previous review works. For example, Chen \textit{et al.}~\cite{chen2020review_UK} and Ji \textit{et al.}~\cite{ji2018review_China} thoroughly examined the recent EVCI developments in the U.K. and China, respectively. LaMonaca's review~\cite{lamonaca2022review} concentrated on introducing stakeholders and policy incentives in the EVCI planning problem. Additionally, studies such as \cite{unterluggauer2022review} highlighted the mutual interests and interconnection between the transport network and the electric grid. The works in \cite{gupta2021review} and \cite{mastoi2022review} detailed various charging technologies and their potential applications for charging station placement. Moreover, Shafiei \textit{et al.}~\cite{shafiei2022review_fastCharging} solely focused on the EVCI planning of fast charging stations. However, due to the expensive upfront costs, fast charging has not been widely spread compared to slow charging. Also, \cite{metais2022review} and \cite{ding2020review} purely analyzed the solution methodologies on EVCI planning strategies. We summarize the key aspects that both previous review papers and our work mainly discuss in Table \ref{tab:comp_previous_review}. In summary, the up-to-date EVCI planning efforts from influential players in the global EV market, e.g., China, the United States, and the European Union, have been inadequately addressed, in particular with the accelerating global EV transition and surging EV uptake. Moreover, reviews on technical planning strategies are not strongly linked to the background knowledge of both stakeholders and charging technologies. However, such background information is essential to establish an efficient EVCI planning model, since the stakeholders' interests are directly related to the objectives of the planning strategies, and the applied charging technologies in EVCI implementation also affect the capacity of charging stations, EV charging speeds, and budgets. Therefore, a systematic review of these aspects is essential. Furthermore, although several studies have delved into the policy incentives and future trends of EVCI planning such as in \cite{gupta2021review,mastoi2022review}, there has been a lack of discussions of specific research gaps and required further contributions on the EVCI planning problem.

\begin{table}[!t]
    \centering
    \caption{The key aspects of EVCI planning discussed in previous review works and our work.}
    \vspace{1mm}
    \resizebox{\columnwidth}{!}{
    \begin{tabular}{c|c|c|c|c|c|c}

    Ref & Year & \makecell{Up-to-date \\ EVCI Practices}  & \makecell{Stakeholders \\ Interests} & \makecell{Charging \\ Technologies} & \makecell{Solution \\ Methodologies} & \makecell{Potential \\ Research Gaps} \\

    \hline 
    
    \cite{chen2020review_UK} & 2020 & \XSolidBrush & \XSolidBrush & \CheckmarkBold & \XSolidBrush & \XSolidBrush \\

    \hline 

    \cite{ji2018review_China} & 2018 & \XSolidBrush & \XSolidBrush & \CheckmarkBold & \XSolidBrush & \XSolidBrush \\

    \hline 

    \cite{lamonaca2022review} & 2022 & \XSolidBrush & \CheckmarkBold & \CheckmarkBold & \XSolidBrush & \XSolidBrush \\

    \hline 

    \cite{unterluggauer2022review} & 2022 & \XSolidBrush & \CheckmarkBold & \XSolidBrush & \CheckmarkBold & \CheckmarkBold \\

    \hline 
    
    \cite{gupta2021review} & 2021 & \XSolidBrush & \XSolidBrush & \CheckmarkBold & \CheckmarkBold & \XSolidBrush \\

    \hline 

    \cite{mastoi2022review} & 2021 & \XSolidBrush & \CheckmarkBold & \CheckmarkBold & \XSolidBrush & \XSolidBrush \\

    \hline 
    
    \cite{shafiei2022review_fastCharging} & 2022 & \XSolidBrush & \XSolidBrush & \CheckmarkBold & \CheckmarkBold & \XSolidBrush \\

    \hline 
    
    \cite{ding2020review} & 2020 & \XSolidBrush & \XSolidBrush & \XSolidBrush & \CheckmarkBold & \XSolidBrush \\

    \hline 

    \cite{metais2022review} & 2022 & \XSolidBrush & \XSolidBrush & \CheckmarkBold & \CheckmarkBold & \XSolidBrush \\

    \hline 

    \textbf{Ours} & 2023 & \CheckmarkBold & \CheckmarkBold & \CheckmarkBold & \CheckmarkBold & \CheckmarkBold 
    
    \end{tabular}
    }
    
    \label{tab:comp_previous_review}
\end{table}

We are motivated to conduct a comprehensive review of the existing efforts on EVCI planning, explore the recent studies of planning strategies, and identify the research gaps that need to be addressed through technical contributions and governmental policy signals in the future. The structure of this paper is organized as follows. Section \ref{sec:practice_and_barriers} revisits the practices of EVCI planning in representative countries in the global EV market, with a further discussion on the potential barriers. Section \ref{sec:stakeholder} offers an overview of the EVCI planning problem by describing the involved stakeholders and their mutual interests. Moreover, Section \ref{sec:charging_tech_and_types} outlines the existing types of charging technologies and their potential implementation in various geospatial locations. We then review solution methodologies of EVCI planning strategies in Section \ref{sec:methodology}. Finally, Section \ref{sec:reserach_gaps} identifies the current research gaps and calls for further contributions to enhance infrastructure planning. We conclude this paper in Section \ref{sec:conclusion}.

\section{Practices and Barriers in EV Charging Infrastructure Planning} \label{sec:practice_and_barriers}

\subsection{Real-world Practices of EVCI Planning}
\label{subsec:practice_and_barriers-practice}

Amid the rapid transition towards EVs worldwide, governments around the world have acknowledged the critical importance of deploying adequate charging infrastructures. They have been proactively funding EVCI development to strengthen and expand their national charging networks, with the aim of meeting the increasing charging demand and promoting sustained EV uptake.

For example, China has emerged as a pioneer in the global EV transition, achieving a remarkable EV market share of $25\%$ of all new vehicles sold nationwide in 2022~\cite{China2023_ICCTreport} To catch the pace of growing EV production, the central government of China has made substantial efforts to support EVCI implementation. By the end of 2021, China's EV charger stock surpassed one million, accounting for $56\%$ of the global total, surpassing Europe and the United States by $2.3$ times and $5.7$ times, respectively~\cite{China2022_ICCTreport}.

Similarly, in the pursuit of net-zero transition, the U.S. Federal Government has set an ambitious goal of reaching $26$ million EVs in stock by 2030 and capturing at least $50\%$ of EV market share across the country by the same year~\cite{us2021_policy}. This requires a significant buildout of public EVCIs, with an expected growth from approximately $0.2$ million in 2020 to $2.4$ million by 2030, signifying an annual increase of $27\%$ charging points~\cite{us2022_NEVI}. For example, California aims to construct better and more equitable EVCI in communities~\cite{us2022_Cali}, which can be achieved by quantifying community EVCI equity values and incorporating them in the designed planning strategies to redistribute charging stations~\cite{Ebbrecht2023,chung2018}.

In addition, the European Union (EU) has urged all member states to build accessible public charging points, enabling EVs to circulate in urban/suburban areas and other densely populated areas~\cite{eu2014_policy}. At the end of June 2023, over $607,000$ public chargers were installed in the EU~\cite{ICCT2023_EUjantojun}, which is still far from its expectation -- $184$ charging points per $100$ kilometers~\cite{eu2016_policyTEN-T}.

\begin{table}[!t]
    \centering
    \caption{The summary of EVCI planning practices.}
    \vspace{1mm}
    \resizebox{\columnwidth}{!}{
    \begin{tabular}{c|c|c|c|c|c}
    
    \makecell{Country/ \\ Organization} & \makecell{EV Market\\ Share} & EVCI Stock & \makecell{Representative \\ EVCI Projects} & Project Target & \makecell{Project \\ Funding} \\
    
    \hline
    \multirow{2}{*}{China} & \multirow{2}{*}{\makecell{$25\%$ by \\ 2022~\cite{China2023_ICCTreport}}} &  \multirow{2}{*}{\makecell{Moret than \\ $1$ million \\ by 2021~\cite{China2023_ICCTreport}}} & \makecell{Pilot EVCI planning \\ in $88$ cities~\cite{China2015_policy}} & \makecell{$1$ charger per $8$ EVs} & \textbackslash \\

    \cline{4-6}
    &&& State Grid program~\cite{China2017_stateGrid} & \makecell{$120,000$ fast chargers \\ $500,000$ public chargers} & \textbackslash \\ 

    \hline

    \multirow{2}{*}{The U.S.} & \multirow{2}{*}{\makecell{$8\%$ by \\ 2022~\cite{IEA2023_EVOutlook}}} & \multirow{2}{*}{\makecell{$0.2$ million \\ by 2021~\cite{ICCT2021_US}}} & \makecell{National Electric \\ Vehicle Infrastructure~\cite{us2022_NEVI}} & \makecell{$500,000$  fast chargers} & \$$7.5$ billion \\

    \cline{4-6}

    &&& \makecell{California Energy \\ Commission program~\cite{us2022_Cali}} & \makecell{Equitable \\ Community EVCIs} & \$$1.4$ billion \\

    \hline

    \multirow{2}{*}{The EU} & \multirow{2}{*}{\makecell{$15.6\%$ by  \\ 2023~\cite{eu2023_marketShare}}} & \multirow{2}{*}{\makecell{$0.6$ million \\ by 2023~\cite{ICCT2023_EUjantojun}}} & \makecell{TEN-T~\cite{eu2016_policyTEN-T} \& CEF-T~\cite{eu2016_policyCEF-T}} &  \makecell{$184$ chargers \\ per 100 km} & \euro$24$ billion \\ 

    \cline{4-6}

    &&& AFIR~\cite{eu2021_AFIR} & \makecell{$1$ kW EVCI  power \\ output per EV} & \textbackslash \\

    \hline

    \multirow{3}{*}{The U.K.} & \multirow{3}{*}{\makecell{$16.6\%$ by \\ 2022~\cite{IEA2023_EVOutlook}}} & \multirow{3}{*}{\makecell{$0.03$ million \\ by 2022~\cite{uk2022_chargerAmount}}} & 2035 Delivery Plan~\cite{uk2020_2035delivery} & \makecell{$6,000$ fast chargers} & \textsterling$950$ million \\ 

    \cline{4-6}

    &&& ORCS~\cite{uk2022_policy} & \makecell{Improve residential \\ charging coverage} & \textsterling$500$ million \\ 

    \cline{4-6}

    &&& \makecell{New building regulation~\cite{uk2021_PM}} & \makecell{Compulsory charger \\ installation} & \textbackslash \\ 

    \hline

    \multirow{2}{*}{Australia} & \multirow{2}{*}{\makecell{$3.8\%$ by \\ 2023~\cite{au2023_nationalStrategy}}} & \multirow{2}{*}{\makecell{$0.005$ million \\ by 2023~\cite{au2023_nationalStrategy}}} & \makecell{National Electric \\ Vehicle Strategy~\cite{au2023_nationalStrategy}} & \makecell{$117$ chargers \\ per 150 km} & \textbackslash \\ 

    \cline{4-6}

    &&& \makecell{VIC~\cite{au2023_nationalStrategy}, WA~\cite{EVC2023}, \\ and SA~\cite{au2022_southAustralia} programs} & \makecell{Statewide EV \\ charging network} & \makecell{AU\$$100$, AU\$$22$, \\ and AU\$$12$ million} \\

    \hline

    \end{tabular}
    }

    \label{tab:country_practice_summary}
\end{table}

We summarize the EVCI planning practices of the above three countries or areas in Table \ref{tab:country_practice_summary}, including their representative policy incentives. The recent progress in the U.K. and Australia is also reviewed, as they are striving to accelerate transport electrification, revealing significant EV market potential.

\subsection{Existing Barriers} \label{subsec:practice_and_barriers-barriers}
Despite the substantial support for EVCI planning and the recent surge in EVCI deployment, the rollout of EVCIs is still significantly lagging behind the pace of EV growth. A report by the Zero Emission Vehicles Transition Council in 2022 emphasized that nearly $6.2$ million public charging points are necessary to accomplish the EV ambition of electrifying the transport sector by 2030~\cite{ZEV2022}. This would correspond to a massive $240$GW of installed EVCI power output. However, as of mid-2022, only $13\%$ of the expected EVCIs were operated, accounting for a mere $10\%$ of the required installed power output. Before pouring more funding or grants into EVCI construction, there are several unresolved problems, which may lead to the stagnant growth of both EVs and EVCIs.

\textbf{Poor EVCI Charging Service:} Home charging via private installed chargers is currently the most preferred option for EV drivers, as they discover that using public EVCIs has evident shortcomings.
\begin{itemize}
    \item \textit{Limited Geographic Coverage and Accessibility:} EV drivers often struggle to find accessible charging stations that align with their daily commuting needs due to the limited geographical coverage of EVCIs. The situation is exacerbated by other uncontrollable factors, such as queuing or infrastructure maintenance.

    \item \textit{Unsatisfactory Charging Experience:} A government report from the U.K. in 2022 highlighted the dissatisfaction of EV drivers with public charging experiences. Common concerns include opaque or excessive charging expenses, poor equipment reliability, and complex access regimes involving numerous apps or smartcards for registration~\cite{uk2022_takingCharge}. Moreover, as few platforms offer transparent and well-rounded information about public EVCIs (such as nearby available chargers, the estimated queuing time, charger's charging speed, and charging prices), EV drivers, with limited information, face difficulties in getting their EVs charged in a reliable, fairly priced, and easily accessed charging station.
\end{itemize}
More importantly, EV drivers without their own driveways heavily rely on public EVCIs and often feel overwhelmed in finding available public chargers, as they have less flexibility in choosing when and where they can charge. In addition, the scarcity of EVCIs deployed in their nearby communities or residential areas forces these EV drivers to accept significantly higher charging fees and endure increased commuting distances compared to those who use private home chargers~\cite{mckinsey2022}. These pain points, combined with the absence of off-street parking and the high upfront costs of private charger installation, may raise consumers' concerns regarding the purchase of EVs.

\textbf{Low Utilization Rates of Public EVCIs:} The reluctance to use public EVCIs, as discussed above, significantly contributes to the low utilization rates of public EVCIs and subsequently disappointing profits from providing charging services. However, ensuring the profitability of public EVCIs (at least enough to cover operation and maintenance costs) is an essential prerequisite for building out a nationwide EV charging network~\cite{mckinsey2022}, as relying solely on funding is not a sustainable long-term solution. The insufficient profitability of deployed chargers can significantly undermine the investment confidence of both governments and utility companies. Consequently, the hesitancy to support further EVCI projects creates a negative feedback loop with EV drivers, thereby failing to deliver more accessible charging points. Such a loop between EVCI users and EVCI developers may ultimately fall into an embarrassing chicken-eggs dilemma, hampering the development of both EVs and EVCIs.

\textbf{Challenging Integration of EVCIs into the Electric Grid:} Efficiently integrating EVCIs into the electric grid poses a significant challenge. While most electric grids can generally meet EV charging demand given the installed generation capacity, few can deliver large amounts of electricity to many EVs or EV fleets at high charging rates concurrently due to the grid's constraints~\cite{eu2022_masterPlan}. This challenge is particularly pronounced in the context of home charging, since the majority of home charging begins after work hours, causing an instant surge in EV charging load and stressing local power distribution networks, which may even lead to severe blackouts. This issue is especially true and complex in countries or regions with high renewable penetration. For instance, in the state of South Australia, rooftop solar photovoltaics can meet all local network demands for more than five hours on a sunny day~\cite{sa2022_rooftopsolar}. Thus, most thermal generators (e.g., coal and gas generators) operate at minimum power output levels or even are shut down. As a result, the grid operator faces significant challenges in matching the rapidly increasing EV charging demand at the onset of home charging. To maintain the grid's security, the excessive charging demand is therefore curtailed. Even though the significant demand can be met, the presence of coal and gas generators pushes up the electricity price and makes home charging less affordable. Such a circumstance may also occur in large-scale concentrated EV fleet charging. Furthermore, as the electric grid is currently transitioning towards low-carbon electricity with increasing renewable uptake, addressing the uncertain surging EV charging loads is not a trivial task given the variability of renewable generation. Upgrading the grid seems to be an option yet extremely expensive. For example, according to the McKinsey market analysis~\cite{mckinsey2022}, for a single fast charging station with four 150 kW fast chargers, the cost of upgrading the grid and site could be more than \$$150,000$.

Apart from the grid security issue associated with increasing EVCI deployment, the upfront cost of connecting the charging stations with the local electric grid may become another hurdle in EVCI implementation, as utility companies, that operate the power distribution network, tend to impose exceptionally high interconnection fees according to the report in~\cite{keith2023}. This might result in EVCI planning that is affected by the considerable interconnection expense, thus leading to location selection sacrificing the convenience and accessibility of public chargers to reduce capital costs. Moreover, the EVCI energization is plagued by notable delays, primarily attributed to complex physical interconnection, difficulties in obtaining easements, and slow permitting processes~\cite{uk2022_takingCharge,IREC2022}. According to the International Renewable Energy Council~\cite{IREC2022}, the average duration for the construction of level 2 and fast chargers can range from one day to six months and from six months to more than two years, respectively. These intensive delays introduce uncertainties for EVCI developers and significantly slow down the transition to electrified transportation.

In addition to necessary governmental support, efficient EVCI planning is inherently essential to breaking down the aforementioned three major barriers for the following reasons.
\begin{itemize}
    \item Efficient EVCI planning involves strategic allocation of available changing resources, effectively enlarging the geographical coverage of built EVCI, resulting in a broader network of charging stations, and providing EV drivers with more charging options. Hence, EV drivers can benefit from accessible and convenient charging services without disruptions to their travel schedules, significantly reducing range anxiety.

    \item By enhancing charging experiences and eliminating range anxiety for EV drivers, efficient EVCI planning can naturally lead to an increase in EVCI utilization rates, which in turn translates to foreseeable profits for EVCI developers. The confidence gained from profitable operations encourages continuous investment in EVCI deployment, creating a positive feedback loop between EVCI developers and EVCI users to break the chicken-and-egg dilemma.

    \item Efficient EVCI planning considers the optimal number of chargers to be connected to the electric grid, alleviating the heavy load burden introduced by large-scale simultaneous EV charging. By ensuring that the surging demand of EV charging does not violate the physical constraints of the grid, efficient infrastructure planning plays a critical role in optimizing the integration of EVs into the existing electric grid.
\end{itemize}
Prior to delving into the review of solution methodologies for devising EVCI planning strategies, establishing a foundational understanding of two pivotal elements, including different stakeholders' interests and existing applicable EV charging technologies, is critical to solving the EVCI planning problem. We therefore introduce and elaborate upon these elements in Sections \ref{sec:stakeholder} and \ref{sec:charging_tech_and_types}, respectively.

\section{Stakeholders and Objectives} \label{sec:stakeholder}
\begin{figure}[!t]
    \centering
    \includegraphics[width=\linewidth]{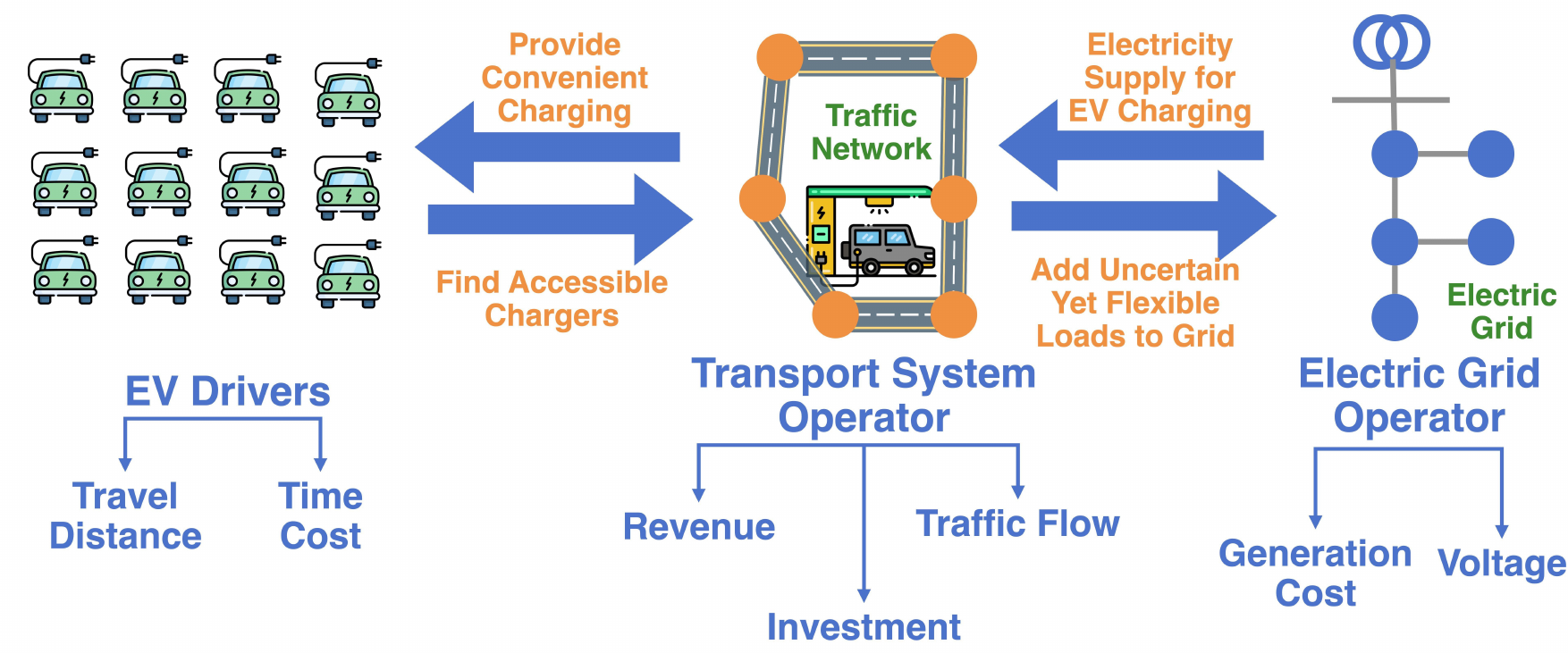}
    \caption{The schematic diagram of EVCI planning with stakeholders.}
    \label{fig:schematic_diagram}
\end{figure}

To overcome the barriers hindering EVCI implementation and devise an efficient solution for EVCI planning, as discussed in Section \ref{sec:practice_and_barriers}, it is crucial to first clarify the key stakeholders involved in the EVCI planning problem, consisting of EV drivers, transport system operators, and electric grid operators. Specifically, the aim of EV drivers is to alleviate the range anxiety problem and get access to convenient charging services, with a focus on optimizing the travel distance and overall time costs throughout the entire charging process. For transport system operators, their primary focus is on improving the utilization rates of EVCIs within budget constraints and generating revenues to cover operation and maintenance expenses. Therefore, their objectives include minimizing investment costs and maximizing captured EV trips and revenue through charging services. As EV charging loads introduce uncertainties to the electric grid, the electric grid operators strive to ensure safe operations (e.g., maintaining stable voltage), while supplying charging electricity with the minimal generation costs. We depict a schematic diagram of EVCI planning with stakeholders in Fig. \ref{fig:schematic_diagram}. Both the interests and objectives of EV drivers, transport system operators, and electric grid operators, are explored in detail in  Section \ref{subsec:stakeholder_EVdriver}, \ref{subsec:stakeholder_transportSystemOperator} and \ref{subsec:stakeholder_electricGridOperator}, respectively.

\subsection{EV Drivers} \label{subsec:stakeholder_EVdriver}
Given that the primary concern of EV drivers centers around the well-known range anxiety problem, research efforts have concentrated on enhancing the accessibility and convenience of charging services. This focus is reflected in the pursuit of the following objectives: 1) minimizing the travel distance to available chargers and 2) reducing overall time costs, consisting of traveling, queuing, and refueling time. Both objectives have been extensively studied in the literature. For example, Hakimi \cite{hakimi1964_pMedian} pioneered an EVCI planning solution aiming to minimize the weighted traveling distance between EV drivers and charging points. This work was subsequently extended in \cite{xu2013_pMedian_BPSO,an2014_pMedian_columnAndConstraint,cavadas2015_MCLM_pMedian_branchAndBound,he2016_pMedian_branchAndBound}. Furthermore, Li \textit{et al.}~\cite{li2023_MCLM_crossEntropy} incorporated the dynamic traffic congestion in the calculation of traveling distances towards charging stations. In addition, Vazifeh \textit{et al.}~\cite{vazifeh2019_SCLM_GA} proposed a data-driven optimization framework to minimize the EV drivers' total excess driving distance to reach charging stations based on real-world traffic data in Boston. Dong \textit{et al.}~\cite{dong2014_MCLM_GA} and Shahraki \textit{et al.}~\cite{shahraki2015_MCLM_branchAndBound} aimed to minimize the number of missed trips due to longer travel distances than the remaining battery range. Also, Wu \textit{et al.}~\cite{wu2023_FCLM_PSO} developed a user-centric EVCI planning strategy to ensure the convenience of EV charging via minimizing detour and travel failure costs. Sadhukhan \textit{et al.}~\cite{sadhukhan2022_onlyGrid} also used the EV power loss (towards charging stations) to describe the accessibility of candidate EVCI planning solutions. Additionally, Pal \textit{et al.}~\cite{pal2023_MCLM} further designed a user convenience factor to characterize the accessibility of charging stations.

Another line of studies draws on minimizing the total time costs of EV charging. The works by Deb \textit{et al.}~\cite{deb2019_MCLM_CSOandTLBO,deb2022_MCLM_CSOandTLBO} developed a multi-objective optimization framework for a coupled traffic-electric networks, where one of the objectives is to reduce the waiting time to get charged. Similar efforts have also been made by Jung \textit{et al.}~\cite{jung2014_MCLM_GA} to minimize the average charging delays in addition to travel time minimization for all EVs on road. Moreover, Pal \textit{et al.}~\cite{pal2023_MCLM} devised a user happiness factor to be maximized, describing the charging congestion, i.e., the queuing length at charging stations. Besides the waiting or queuing time of EV charging, several studies have taken the EV traveling time and recharging time into consideration. For instance, He \textit{et al.}~\cite{he2015_MCLM_GA} positioned charging stations within the driving ranges of EVs (based on simulated battery's state of charge), while minimizing the overall travel and refueling time. Studies, such as~\cite{ma2021_FCLM_branchAndBound,li2023_MCLM_crossEntropy}, provided a more detailed perspective of the time costs, covering traveling, queuing, and recharging time as three sub-objectives to be optimized. As the queuing time directly impacts the charging experiences of EV drivers, Zhao \textit{et al.}~\cite{zhao2021_onlySizing_GA} located charging stations to improve the quality of service of EVCIs via fuzzy decision-making.

\subsection{Transport System Operator} \label{subsec:stakeholder_transportSystemOperator}
Transport system operators often serve as the EVCI developers and are often funded by governments for the deployment of charging stations. Therefore, the primary focus is typically on minimizing the infrastructure costs associated with EVCI planning, encompassing investment, operation, and maintenance costs. Numerous studies have delved into minimizing such infrastructure costs, as seen in \cite{jia2012_SCLM_branchAndBound,liu2013_onlyGrid,sadeghibarzani2014_SCLM_GA,xiang2016_FCLM_branchAndBound}. 

In addition, driven by the increasing adoption of EVs and the imperative to improve the utilization rates of EVCIs, Toregas \cite{toregas1971_SCLM} introduced a classic facility location model with the context of traffic networks, aiming to cover all charging demand in the network. Variants of this method have been widely explored in studies such as \cite{church1974_MCLM,xi2013_MCLM_branchAndBound,cai2014_MCLM,asamer2016_MCLM_branchAndBound}, where the number of installed charging stations is constrained by budget limitations. Previous approaches tend to analyze the static view of EV drivers' charging demand at traffic nodes, which is often represented as varying traffic flows, i.e., the origin-destination EV trips, instead. To derive better EVCI planning solutions, a line of work has focused on maximizing captured traffic flows considering budget constraints~\cite{hodgson1990_FCLM,kuby2005_FRLM,kim2012_FRLM_branchAndBound} and capacity limits of each charging station~\cite{upchurch2009_CFRLM,li2016_CFRLM_GA,wang2018_CFRLM_GA}.

Given the substantial investment costs of EVCI implementation, which have been increasingly covered by third-party EVCI developers rather transport system operators, these developers have developed smart charging strategies or flexible charging scheduling to maximize operational revenues through charging services or vehicle-to-grid technologies, capitalizing on fluctuating electricity prices~\cite{chen2019_privateStationInvestorPerspective,zhao2020_privateStationInvestorPerspective}. For example, the Boston Consulting Group explored strategic approaches to deploying charging infrastructures, including the involvement of third-party developers and equipment suppliers~\cite{bcg2021}. Though EVCI developers can, to some extent, accommodate the increasing EV charging demand, their revenue-driven planning strategies may lead to imbalanced distributions of charging stations. For example, these developers tend to invest more in areas where individuals can afford their designed charging pricing, e.g., high-income communities. In response to these market-driven initiatives, the transport system operators are increasingly engaging in public-private partnerships with these developers for sustainable EVCI construction~\cite{mckinsey2022,zhang2018_ppp,su2016_ppp}. In this context, the role of the transport system operator evolves towards enhancing social welfare, focusing on equitable EVCI placement and accessibility.

\subsection{Electric Grid Operator} \label{subsec:stakeholder_electricGridOperator}
Given the strong coupling relationship between the traffic network and the electric grid, the increasing EV penetration in the transport sector inevitably amplifies the load in the coupled power distribution network, leading to voltage deviations and posing threats to the reliability of electricity supply~\cite{davidov2017_SCLM_branchAndBound}. Though uncertainties associated with EV charging, such as stochastic EV arrivals and variable charging durations, can be transformed into beneficial flexibilities with proper management, as supported by Madzharov \textit{et al.}~\cite{madzharov2014}. At the planning stage, the potential impacts of EV charging are commonly considered as an unresolved uncertainty that may affect the grid's ability, e.g., in worst-case scenario analysis, as EV charging behavior is uncontrollable. Therefore, satisfying the operational constraints of the distribution network becomes a critical prerequisite for EVCI planning in coupled networks from the perspective of electric grid operators. Correspondingly, the works by \cite{zhang2018_CFRLM_branchAndBound,xiao2023_FCLM_GA} incorporated distribution network constraints, e.g., the active/reactive power flow balance, voltage, and line capacity bounds. Lei \textit{et al.}~\cite{lei2022_onlySizing} further employed an alternative-current power flow model (with the objective of minimizing the total generation costs) to model the dynamics of the distribution network. More studies tend to focus on maintaining the grid's voltage stability, with works such as \cite{wang2013_FCLM_crossEntropy,mao2021_FCLM_crossEntropy,deb2022_MCLM_CSOandTLBO} minimizing the voltage deviations across the distribution network. Optimizing the overall energy losses has also been carried out as part of the EVCI's investment costs in studies such as \cite{liu2013_onlyGrid,yao2014_FCLM_MOEAD,deb2022_MCLM_CSOandTLBO}. While previous methods require full knowledge of the distribution network parameters, Li \textit{et al.}~\cite{li2023_MCLM_crossEntropy} addressed uncertain load burdens (caused by EV charging) in a realistic region-based distribution network, without the need for detailed topology and parameter information. Despite introducing additional uncertainty to the electric grid, Das \textit{et al.}~\cite{Das2020_chargingtype} demonstrated that effectively scheduled EV charging is able to reduce energy peaks of the power systems, i.e., peak shaving, and even enhance grid resilience without altering the underlying structure of the distribution network. However, reinforcing the grid may still be essential to accommodate the growing number of EVs~\cite{Sridhar2022_chargingtype}. Sridhar \textit{et al.} \cite{Sridhar2022_chargingtype} studied the optimal number of transformers to upgrade every five years based on different rates of investment. The research also indicated that upgrading transformers may not be necessary if charging is carefully managed by the electric grid operator and transformer capacities are not breached, since such grid reinforcement would entail extra investment in communication technologies~\cite{Brinkel2020_chargingtype}.

\section{Charging Technologies and Types} \label{sec:charging_tech_and_types}

In addition to stakeholders' interest, various charging technologies are also evolving with the rapid global EV transition. This section aims to provide a detailed description of existing charging technologies and their potential applications in real-world EVCI planning problems. The structure of this section is organized as follows. Section \ref{subsec:charging_tech_and_types_chargingLevel} delves into the technical aspects of level 1, 2, and 3 of EV charging, including corresponding charging modes and types, followed by Section \ref{subsec:charging_tech_and_types_chargingLocation}, exploring the applications of various charging technologies in the literature. Finally, in Section \ref{subsec:charging_tech_and_types_matchTypeAndTech}, we match charging technologies to potentially appropriate areas, such as residential areas and workplaces, for effective EVCI rollout.

\subsection{Level 1, 2, and 3 Charging Technologies and Features}
\label{subsec:charging_tech_and_types_chargingLevel}
\subsubsection{Charging Level}
Charging levels 1, 2, and 3 refer to the speeds of EV charging, with each level offering varying rates of charging. Level 1, characterized by the slowest charging speed, is commonly employed at household outlets, enabling EV drivers to incrementally recharge their batteries~\cite{EVC2023_chargeType}. Charging level 2, comparatively faster than level 1, requires additional installation. At residences, it operates at power rates of $7.4$ kW (i.e., single-phase chargers), while at public charging stations, it ranges from $11$ to $22$ kW (i.e., three-phase chargers). The fastest charging speed is achieved at fast charging stations, denoted as level 3 charging, with a rated power of up to $350$ kW. According to a report by the transport and road agency of the New South Wales government, namely Transport for NSW, one hour of level 3 charging can support approximately $150$-$300$ kilometers driving range~\cite{NSW2023}. 

\begin{table}[!t]
    \centering
    \caption{The characteristics of level 1, 2, and 3 EV charging technologies.}
    \resizebox{\columnwidth}{!}{
    \begin{tabular}{c|c|c|c|c|c|c}    
        Level & \makecell{AC/DC \\ Charging} & \makecell{Rated \\ Power} & \makecell{Charging \\ Time} & Charging Rate & Mode & Type \\
        
        \hline 
        
        Level 1 & AC &  $1.4$-$3.7$ kW & $5$-$16$ hours & $10$-$20$ km/hour & Mode 2 & Type 1 \\
        
        \hline 
        
        \multirow{2}{*}{Level 2} & \multirow{2}{*}{AC} & \makecell{$7.4$ kW \\ single-phase}  & $2$-$5$ hours & $30$-$45$ km/hour & \multirow{2}{*}{Mode 3} & Type 1 \\

        \cline{3-5} \cline{7-7}
        
        && \makecell{$11$-$22$ kW \\ three phase} & \makecell{$30$ mins \\ to $2$ hours} & $50$-$130$ km/hour & & Type 2 \\
        
        \hline 
        
        Level 3 & DC & $25$-$350$ kW & $10$-$60$ mins & $150$-$300$km/hour & Mode 4 & \makecell{CHAdeMO \\ CCS2} \\
        
    \end{tabular}
    \label{tab:charging_levels}
    }
\end{table}

\subsubsection{Charging Mode}
Charging mode is defined as the safety communication protocol between the EV and the charger \cite{EVHub2020_chargingMode}, falling into three types. In Mode 1, an EV is directly connected to the household socket, and there is no communication or data exchange between the EV and the charger. In contrast, Mode 2 improves safety functions during the charging process with additional in-cable control and safety devices. Safety features include limiting the current flow without the verification of a protective earth connection to ensure that the EV has been safely plugged into the charger, thereby preventing over-current or over-temperature connections~\cite{thomas2012_chargingMode}. Mode 3 has a similar communication protocol to Mode 2, while enabling higher power transfer. This mode can charge between $3.6$ kW (in single-phase chargers) to $40$ kW (in three-phase chargers), determined by various types of cables and charging capacity of EVs. Different from the previous three modes, Mode 4 is designed for direct current (DC) charging, i.e., fast charging. In this mode, the DC charger is integrated inside the wall box of the charging station, allowing the charging current to directly flow into the EV's battery. In the context of alternative current (AC) charging, the current flows into the EV's charger rather than directly into the battery. Therefore, the safety protocols of Mode 4 are more strict than the aforementioned modes. 

\subsubsection{Charging Type}
Charging type refers to the specification of the socket outlet~\cite{FCAI2019}. For AC charging, both Type 1 and Type 2 charging plugs can be utilized. Type 1 is a single-phase plug permanently connected to the charging station, and enables up to $7.4$ kW of charging power. Type 2 has a three-phase plug and supports higher power rates of up to $22$ kW at home or $43$ kW at public charging stations. With regards to DC charging, two representative plugs are widely adopted, referred to as the CHAdeMO and the CCS2. The former, which was developed in Japan, allows for both fast charging and bi-directional charging with a maximum charging power of $100$ kW. The CCS, an extension of the Type 2 plug, supports both AC and DC charging within $350$ kW. Charging types are adopted by countries based on their national EV uptake. For example, Australia uses Type 2 for AC charging and provides adaptors (from Type 1 to Type 2) for older models using legacy Type 1 plugs. For DC charging, although both CHAdeMO and CCS2 are viable, CCS2 is preferred by EV drivers due to its compatibility with Type 2 plugs. In comparison, CHAdeMO is only compatible with a few Japan-manufactured models.

We summarize the characteristics of the three levels of EV charging, including their rated power, expected charging durations, charging rates, charging modes, and their corresponding charging types, as shown in Table \ref{tab:charging_levels}.

\subsection{Applications of Charging Technologies in EVCI Planning}
\label{subsec:charging_tech_and_types_chargingLocation}
\subsubsection{Planning of Single Type of Charger}
Numerous studies in the literature have focused on singular types of chargers, e.g., level 2 or 3 charging. For instance, Babic \textit{et al.}~\cite{Babic2022_chargingtype} optimized the number of parking spots to be converted into charging spots equipped with level 2 chargers in a given parking lot, which was then generalized for semi-residential areas, shopping centers, airports, and homogeneous locations where EV are parked for an extended period. Level 2 charging was also considered in the context of a university campus~\cite{Abdullah2022_chargingtype}, where EVs typically remain parked during working hours. Regardless of the significantly shorter charging time provided by fast charging, level 2 charging is preferred due to the common parking behaviors of most EV drivers. Sathaye \textit{et al.}~\cite{Sathaye2013_chargingtype} employed separate methods for EVCI planning of level 2 and level 3 charging stations on highway corridors. For level 2 charging, the optimal density of charging stations was determined based on population density within a preset budget. In contrast, they optimized the density of level 3 chargers by adjusting the input parameter of geospatial distance to the nearest charging ports, which was not subject to budget limitation.

While level 2 chargers require lower investment than level 3 charging, the latter has drawn increasing attention due to their more effective and efficient service capability for large amounts of EVs~\cite{Nie2013_chargingtype}. Given the scarce distribution of chargers in the current traffic network, battery capacity consumed while driving to charging stations (termed \textit{EV loss} in~\cite{sadeghibarzani2014_SCLM_GA}) and the considerable infrastructure costs are two major drivers hindering the widespread of EVCI planning. To address these barriers, studies such as~\cite{Kong2019_chargingtype,Chen2020_chargingtype} have explored the planning of fast chargers on roads to mitigate range anxiety with better quality of service.

\subsubsection{Planning of Multi-Type of Chargers}
Despite the in-depth research in EVCI planning with a single type of charger, emerging trends favor the incorporation of multiple types of chargers in various geographical locations, e.g., residential, commercial, and industrial areas, for a more holistic EVCI planning solution~\cite{Patil2023_chargingtype}. Given that EV drivers often commute between multiple areas, the deployment of different types of chargers can complement each other, providing more effective service for EV charging demand~\cite{xi2013_MCLM_branchAndBound}. For example, studies in~\cite{Luo2018_chargingtype,Rene2023_chargingtype,Jiang2023_chargingtype} considered EVCI planning in both residential and commercial areas. In these studies, level 1 and 2 chargers are utilized in residential areas, while a mix of level 2 and DC charging are designed for commercial areas for fast charging speeds. Liu \textit{et al.}~\cite{Liu2022_chargingtype} extended such planning with additional inclusion of industrial areas. Moreover, Zhang \textit{et al.}~\cite{zhang2016_SCLM_PSO} considered EVCI planning in a public parking lot with level 2 charging, along with roadside level 3 charging within an urban traffic network, aiming to find the optimal number of placed chargers of the network. Similarly, Wang \textit{et al.}~\cite{wang2013_FCLM_branchAndBound} introduced slow charging (i.e., level 1 and 2 charging) for attractions, fast charging at roadsides, and mixed chargers at convenience stores in a real-world traffic network. 

Beyond level 2 and 3 chargers, He \textit{et al.}~\cite{he2015_MCLM_GA} included level 1 chargers in their user-centric EVCI planning, capturing EV driving behaviors. Specifically, EV drivers were simulated to stop at different locations. The placement of charging stations was then optimized to mitigate the inconveniences of charging services, referred to as social welfare costs. Following the activity-based simulations, Dong \textit{et al.}~\cite{dong2014_MCLM_GA} assumed that all drivers have access to level 1 charging at home and can choose to charge either at workplaces or roadsides with level 2 charging. Moreover, Schoenberg \textit{et al.}~\cite{schoenberg2023_onlySizing} differentiated slow and fast charging for destination and en-route charging, respectively, with the objective of minimizing EV waiting times. With considerations of multiple stakeholders, Li \textit{et al.}~\cite{li2021_privateStationInvestorPerspective} determined the optimal number of level 2 and 3 chargers with varying budgets under scenarios of different mixes of chargers, from the perspectives of both EV drivers and private EVCI investors. This study suggested that, for investors with lower budgets, level 2 charging comprises the majority of charging stations in the traffic network, with few stations dedicated to level 3 charging. If budget permits, investors prefer building level 3 charging stations given higher profit margins. However, the research by Madina \textit{et al.}~\cite{Madina2016_chargingtype} indicated that a high utilization rate is an essential prerequisite for the implementation of level 3 charging, ensuring a positive return on investment, given the high upfront costs.

In addition to EVCI planning targeting private EVs, few studies have been conducted on building chargers for public EVs, e.g., electric taxis~\cite{cai2014_MCLM,jung2014_MCLM_GA,asamer2016_MCLM_branchAndBound} and electric buses~\cite{wang2022_onlySizing}. The works by \cite{cai2014_MCLM} and \cite{jung2014_MCLM_GA} aimed to design an optimal mix of level 1, 2, and 3 chargers for large-scale electric taxi fleets, while Asamer \textit{et al.}~\cite{asamer2016_MCLM_branchAndBound} focused on fast charging. Additionally, to mitigate the urban-rural divergence, EVCI planning with suitable types of chargers has also been discussed in the literature. According to the work by McKinney \textit{et al.}~\cite{mckinsey2022}, a rural neighborhood can experience the minimum impacts on its fragile connected grid if all vehicles are electrified and also have access to level 2 home charging. To achieve that, Baidiei \textit{et al.}~\cite{Badiei2023_chargingtype} proposed a community-based approach, with economic-feasible guarantees, to upgrade vehicles and incorporate EVs into rural electric grids. Specifically, customers can own or subscribe to solar arrays that form solar-powered charging stations and receive credits for excess solar generated. Utilizing solar energy to power charging stations in rural areas is highly suitable as these charging stations generally require large amounts of space and thus, are not deemed suitable for urban locations~\cite{Das2020_chargingtype}.

\subsection{Match Charging Technologies to Various Geographical Areas}
\label{subsec:charging_tech_and_types_matchTypeAndTech}
As discussed above, there are clear trends regarding the suitability of specific charger types at different locations. In this section, we categorize charging location types into residential (in Section \ref{subsubsec:charging_tech_and_types_matchTypeAndTech_residentialCharging}), destination (in Section \ref{subsubsec:charging_tech_and_types_matchTypeAndTech_destinationCharging}), and public (in Section \ref{subsubsec:charging_tech_and_types_matchTypeAndTech_enRouteCharging}) EV charging, with detailed discussions on relevant charging technologies. 

\subsubsection{Residential Charging}
\label{subsubsec:charging_tech_and_types_matchTypeAndTech_residentialCharging}
Level 1 charging is the most prevalent charging type at residences and is preferred over faster charging options, i.e., level 2 and 3 charging, due to its affordability \cite{Madina2016_chargingtype}. As drivers naturally stay parked in residential areas for extended periods, often overnight, fast charging seems to be unnecessary compared to level 1 charging. By plugging in their EVs at home, EV drivers can replenish the battery used during their daily, local journeys. Studies that integrate user journeys, such as \cite{he2015_MCLM_GA} and \cite{dong2014_MCLM_GA}, often assumed that drivers have access to level 1 residential charging. Although level 1 charging is sufficient for urban residential areas, McKinney \textit{et al.} \cite{McKinney2023_chargingtype} explored the potential of utilizing level 2 charging in rural homes. Even with higher charging rates, level 2 charging can be managed with minimal grid impact in the case that all rural households have access to level 2 home charging.

\subsubsection{Destination Charging}
\label{subsubsec:charging_tech_and_types_matchTypeAndTech_destinationCharging}
Destination EV charging includes areas where drivers finish their journey or stop for an extended period of time, such as at a workplace, university campus, shopping center, and roadsides near their destinations. Level 2 charging has been mostly discussed in the context of these locations. EVs often take approximately $2$ to $5$ hours for charging completion, which is less than or at most matches the duration that drivers usually stay at these locations. These stations are often integrated into existing parking structures and locations with longer dwell times, necessitating coordination with municipal planning efforts. As the charging rate offered by level 2 charging differs compared to level 3 charging, urban charging demand, therefore, has a relatively moderate impact on the local distribution network, except during peak urban demand hours~\cite{wu2022fail}, e.g., surge of workplace charging in the morning. However, Nie \textit{et al.}~\cite{Nie2013_chargingtype} revealed that level 2 charging is not the optimal solution for the surging EV uptake, while level 3 charging is preferred in future scenarios if its construction is economically feasible. A viable solution that caters to a diverse driver group is to build charging stations combining AC and DC charging~\cite{Liu2019_chargingtype}, which can potentially reduce the overall number of charging stations needed to be installed, as well as the total investment costs.

\subsubsection{En-Route Charging}
\label{subsubsec:charging_tech_and_types_matchTypeAndTech_enRouteCharging}
In contrast to residential and destination EV charging, EV drivers inevitably stop at charging stations during their trips due to the limitation of EV driving ranges, referred to as en-route charging. Being able to charge quickly, such that EV drivers can finish their trips with refueled EVs without any disturbances on their traveling plans. In the literature, studies predominantly focused on implementing level 3 chargers on roadsides~\cite{sadeghibarzani2014_SCLM_GA,Kong2019_chargingtype,Chen2020_chargingtype} or highways~\cite{Sathaye2013_chargingtype,Madina2016_chargingtype,Micari2017_chargingtype}, given the ability of fast charging to charge EVs up to $80\%$ within one hour, satisfying the urgent charging needs of EV drivers. Such deployment often requires integration with transportation planning for long-distance routes and corridors, often in conjunction with rest stops and service areas. Different from level 2 charging, level 3 charging, due to its significant charging power, has a considerable adverse impact on the electric grid, particularly during high usage periods. Though building charging stations (equipped with fast chargers) requires substantial infrastructure investment costs, the optimal number of charging stations to be installed can be reduced, as fast charging can provide longer battery ranges within shorter charging time compared to other levels of chargers.

In summary, level 1 charging is commonly employed in residential areas and private households. Level 2 charging is a more suitable choice for workplaces, university campuses, or residential parking lots, where EVs tend to park for extended periods, taking up to five hours to get fully charged. Level 3 is often located on roadsides or highways to meet the en-route charging demand, significantly reducing expected charging time.

Table \ref{tab:charging_locationtypes} matches charging levels to related location types.

\begin{table}[!t]
    \centering
    \caption{Suitable location types for respective charging levels.}
    \begin{tabular}{c|c|c|c}
        Ref. & Charger & Rated Power & Planning Location \\
        
        \hline
        
         \cite{Abdullah2022_chargingtype} & Level 2 & $6.6$ kW & Campus parking \\

         \hline
         
         \cite{Babic2022_chargingtype} & Level 2 & $7.7$ kW & Semi-residential parking \\
         
         \hline
         
         \cite{McKinney2023_chargingtype} & Level 2 & $7$ kW & Rural homes\\
         
         \hline
         
         \multirow{2}{*}{\cite{shahraki2015_MCLM_branchAndBound}} & Level 2 & $7$ kW & \multirow{2}{*}{Roadside} \\ 

         \cline{2-3}
         
         & Level 3 & $37.5$ kW \\

         \hline
         
         \multirow{2}{*}{\cite{Liu2019_chargingtype}} & Level 2 & $7$ kW & \multirow{2}{*}{Roadside} \\ 
         
         \cline{2-3}
         
         & Level 3 & $45$ kW \\
         
         \hline
         
         \multirow{2}{*}{\cite{Micari2017_chargingtype}} & Level 2 & $22$ kW & \multirow{2}{*}{Highway} \\ 

        \cline{2-3}
         
         & Level 3 & $50$ kW \\
         
         \hline
         
         \cite{sadeghibarzani2014_SCLM_GA} & Level 3 & $50$-$250$ kW & Roadside \\
         
         \hline
         
         \cite{Kong2019_chargingtype} & Level 3 & $40$ kW & Roadside \\
         
         \hline
         
        \cite{Madina2016_chargingtype} & Level 3 & $50$ kW & Highway \\
        
         \hline
         
         \cite{wang2009_SCLM_branchAndBound} & Level 3 & Not specified & Highway \\
         
         \hline
         
         \cite{Dong2016_chargingtype} & Level 3 & Not specified & Highway \\
         
         \hline
         
         \cite{Wang2018_chargingtype} & Level 3 & Not specified & Highway \\
         
         \hline
         
         \cite{Chen2020_chargingtype} & Level 3 & $90$ kW & Roadside \\
         
    \end{tabular}
    \label{tab:charging_locationtypes}
\end{table}

\section{Solution Methodology for Planning} \label{sec:methodology}

EVCI deployment has been widely treated as a facility location problem, categorized into two primary solution methods: 1) node-based approaches and 2) flow-based approaches. With traffic networks represented as graphs consisting of nodes and edges (or links), node-based approaches involve siting charging stations at traffic nodes, functioning as centroids to fulfill EV charging demands radiating outward while minimizing the overall infrastructure costs. On the other hand, flow-based approaches position EVCIs along the EV's origin-destination (O-D) trips, targeting maximum captured EV traffic flows. Multifaceted challenges in the EVCI planning problem further complicate the facility location problem, including fluctuating EV charging demands over time, substantial infrastructure costs, potential adverse impacts on the electric grids due to the inherent interconnection between traffic networks and the electric grid, as well as the intricate task of determining the optimal capacity for charging stations with specific chargers or a mix of various types of chargers. The structure of this section is as follows. In Section \ref{subsec:methodology_preliminaryModel}, we introduce fundamental models of EVCI planning in traffic networks, including both node-based and flow-based models investigated in Section \ref{subsubsec:methodology_preliminaryModel_nodeBasedModels} and \ref{subsubsec:methodology_preliminaryModel_flowBasedModels}, respectively. Addressing key challenges, such as capacity sizing of charging stations and strategies for modeling the uncertain routing choices of EV traffic flows, is then detailed in Section \ref{subsubsec:methodology_preliminaryModel_capacitySizing} and \ref{subsubsec:methodology_preliminaryModel_uncertainEVFlowModeling}, respectively. Furthermore, recognizing the inherent interdependency between traffic networks and the electric grid, we delve into the existing literature on strategies for EVCI planning in the coupled systems in Section \ref{subsec:methodology_strategyInCoupledNetworks}.

\subsection{Preliminaries of EVCI Planning} \label{subsec:methodology_preliminaryModel}

\subsubsection{Node-based Models} \label{subsubsec:methodology_preliminaryModel_nodeBasedModels}

Toregas~\cite{toregas1971_SCLM} introduced the set covering location model (SCLM) with the aim of minimizing the total infrastructure costs for EVCI deployment while ensuring coverage for all charging demands. This approach concentrates charging demands at specific traffic nodes, strategically placing EVCIs to guarantee that no demand node exceeds a predetermined geospatial service distance. Due to budget constraints in EVCI planning, Church~\cite{church1974_MCLM} extended the SCLM into the maximum covering location model (MCLM). The MCLM, similar to the SCLM in considering service distance, allowed demand nodes not to be covered for budget limitations. Different from defining service distance for charging coverage, Hakimi~\cite{hakimi1964_pMedian} proposed the p-median model to place a total of $p$ charging points to minimize the weighted traveling distance between demand and charging points.

In literature, Wang \textit{et al.}~\cite{wang2009_SCLM_branchAndBound} employed the SCLM model to place fast charging stations for intercity EV trips in Taiwan, which was then extended into a combination of SCLM and MCLM models considering intra-city EV charging in~\cite{wang2010_SCLM_branchAndBound} and solved by the branch-and-bound (B\&B) algorithm. Building on these studies, similar investigations were conducted in the traffic networks of Beijing~\cite{liu2012_SCLM_APSO}, Stockholm~\cite{jia2012_SCLM_branchAndBound}, and Boston~\cite{vazifeh2019_SCLM_GA}. Also, Wang \textit{et al.}~\cite{wang2017_SCLM} and Davidov \textit{et al.}~\cite{davidov2017_SCLM_branchAndBound} demonstrated the feasibility of the SCLM model in large-scale synthetic traffic networks. Additionally, various studies delved into the MCLM \cite{xi2013_MCLM_branchAndBound,cai2014_MCLM,dong2014_MCLM_GA,jung2014_MCLM_GA,he2015_MCLM_GA,shahraki2015_MCLM_branchAndBound,tu2016_MCLM_GA,asamer2016_MCLM_branchAndBound} and p-median \cite{xu2013_pMedian_BPSO,an2014_pMedian_columnAndConstraint,he2016_pMedian_branchAndBound} models, some of which were tackled by heuristic algorithms, e.g., particle swarm optimization (PSO) algorithm~\cite{xu2013_pMedian_BPSO} and genetic algorithm (GA)~\cite{dong2014_MCLM_GA,jung2014_MCLM_GA,he2015_MCLM_GA,tu2016_MCLM_GA}, due to their faster computation speed in such an NP-hard problem. Furthermore, He \textit{et al.}~\cite{he2016_pMedian_branchAndBound} conducted a comparative analysis of the SCLM, MCLM, and p-median models, revealing that the p-median model excelled in producing efficient planning solutions, as it tended to position charging stations in close proximity to communities with higher EV demand.

A significant limitation of node-based models is their static vision of charging services. SCLM and MCLM operate under the assumption that the service distance of each charging station remains constant. This assumption implies that a charging station can only cater to demand points within its fixed service radius, which does not hold in practice and is not always accurate, since it is influenced by the unique, uncertain, and uncontrollable charging behaviors of EV drivers. The p-median model, while not reliant on predefined service distances, also falls short by not incorporating realistic EV charging data. Consequently, all these node-based models tend to result in static and less efficient EVCI planning.

\subsubsection{Flow-based Models} \label{subsubsec:methodology_preliminaryModel_flowBasedModels}
In node-based models, besides simplifying the variable service distance of charging stations, these models also rely on estimating charging demands and assuming that such demands are aggregated at traffic nodes. However, in real-world scenarios, EV charging demands are not always explicitly defined at nodes; instead, they are often represented as traffic flows in the form of O-D pairs~\cite{hodgson1990_FCLM}. To address the limitations of node-based models, Hodgson~\cite{hodgson1990_FCLM} developed a flow-based adaptation of the MCLM, referred to as the flow-capturing location model (FCLM). The FCLM formulates O-D trips of EVs and aims to maximize captured EV flows on the shortest path between origins and destinations. In FCLM, one O-D pair is considered covered if it passes through at least one traffic node equipped with a charging station. Since an EV may need to stop at multiple charging stations to complete its O-D trip, Kuby~\cite{kuby2005_FRLM} extended the FCLM into the flow-refueling location model (FRLM) to place a set (or combination) of charging stations along O-D trips to meet such refueling needs. Research by Upchurch \textit{et al.}~\cite{upchurch2010_nodeBased} demonstrated that flow-based methods tend to be more stable in EVCI deployment as the number of charging stations to be installed increases. Case studies utilizing FCLM and FRLM models are presented in \cite{wang2013_FCLM_branchAndBound,huang2015_FCLM_branchAndBound,sun2020_FCLMfast_MCLMslow_branchAndBound} and \cite{kim2012_FRLM_branchAndBound,wu2017_FRLM_branchAndBound}, respectively.

\subsubsection{Capacity Sizing of Charging Stations} \label{subsubsec:methodology_preliminaryModel_capacitySizing}

Notably, both flow-based and node-based models typically assume that installed charging stations have unlimited capacity to serve incoming EVs, which may not align with practical constraints, as budget considerations and actual charging demands often limit the number of chargers at a single station. Moreover, a single combination of charging stations in the FRLM model (or a single charging station in the FCLM model) may fail to refuel all traffic flows of a particular O-D pair when introducing the capacity constraint. 
To address this issue, Upchurch~\cite{upchurch2009_CFRLM} improved the FRLM by incorporating capacity sizing for candidate charging sites, namely the capacitated FRLM (CFRLM). Specifically, in the formulation of CFRLM, the binary variable at each traffic node (i.e., determining the placement of charging stations) is replaced with an integer variable indicating the number of chargers to be installed at the specific candidate sites, while the total number of chargers (that can be implemented) is subject to budget considerations. The CFRLM avoids oversimplifying the implementation of EVCIs and therefore has been leveraged in multiple studies for more practical EVCI planning strategies~\cite{li2016_CFRLM_GA,wang2018_CFRLM_GA,zhang2018_CFRLM_secondOrderCone}

The introduction of capacity sizing enhances the applicability of planning solutions in real-world EVCI deployment scenarios. Similar improvements have also been integrated into node-based models. For example, Dong \textit{et al.}~\cite{dong2016_SCLM_SNN} assumed the number of chargers of each candidate charging station should satisfy the peak number of EVs to be charged in their SCLM model. Moreover, Rajabi \textit{et al.} and Cavadas \textit{et al.} constrained the number of deployed fast~\cite{rajabi2017_SCLM_GA} and slow~\cite{cavadas2015_MCLM_pMedian_branchAndBound} chargers in urban areas, respectively, while Zhang \textit{et al.}~\cite{zhang2016_SCLM_PSO} considered different types of capacitated charging facilities, including level 1, 2, and 3 chargers, to minimize the overall social costs of the whole EV charging system. Furthermore, several studies focused on the capacity sizing problem of existing EVCI, e.g., Zhao \textit{et al.}~\cite{zhao2021_onlySizing_GA} solved the sizing problem by maximizing the fuzzy quality of service of charging stations via the alpha-cuts-based algorithm, Lei \textit{et al.}~\cite{lei2022_onlySizing} determined how many charging units should be installed at different locations on the intercity traffic network from the perspective of carbon emission, Wang \textit{et al.}~\cite{wang2022_onlySizing} investigated the optimal capacity of two types of charging points for electric buses (i.e., destination charging and en-route charging), Schoenberg \textit{et al.}~\cite{schoenberg2023_onlySizing} discussed the optimal sizing of existing EVCI to potentially improve the average extra time spent with charging for EV drivers on everyday trips, and Duan \textit{et al.}~\cite{duan2022_onlySizing} discussed the fast charger sizing at specific locations considering the market competition from other charging service providers.

In addition to directly incorporating capacity sizing into the mathematical formulation of EVCI planning, recent studies have drawn increasing attention to the queuing theory to address this specific sizing problem. Charging facilities have their capacities that can serve a certain number of EVs, and they cannot refuel EVs instantaneously, thereby resulting in charging congestion and EV queuing. To characterize the queuing behaviors of EVs, recent studies~\cite{xiang2016_FCLM_branchAndBound,yang2017_MCLM_branchAndBound,yang2019_onlySizing_MarkovDemand} modeled the charging station as a multiple-service queuing system, while the arrival of EVs was described as a Poisson process and often estimated using historical EV travel data. The capacity determination is then converted into the objective of minimizing the queuing time of EV charging~\cite{li2023_MCLM_crossEntropy}, which is co-optimized with the initial objectives of minimizing infrastructure costs or maximizing the captured EV traffic flows.

\begin{figure}[!t]
    \centering
    \includegraphics[width=\linewidth]{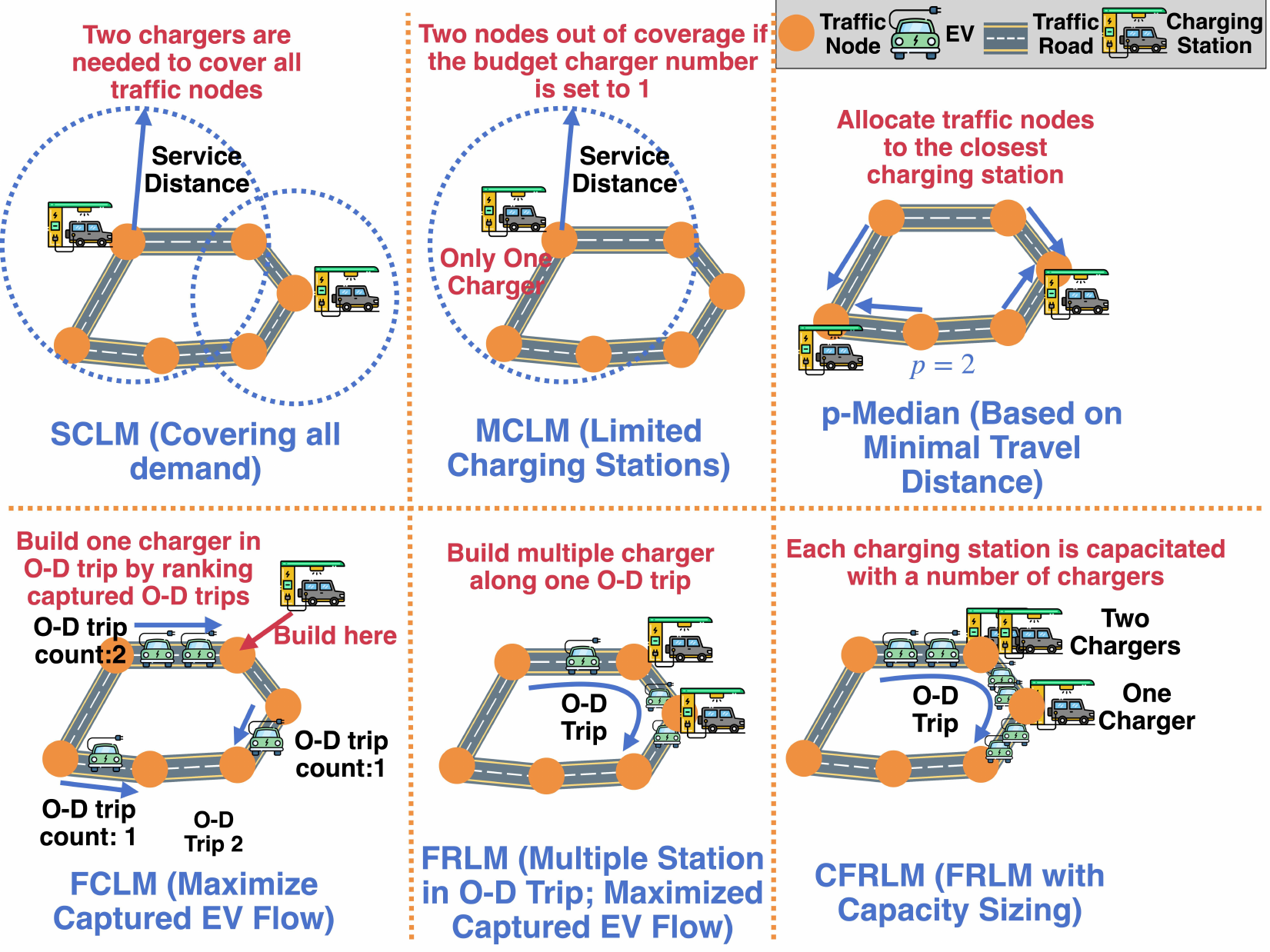}
    \caption{The illustrations of SCLM, MCLM, p-median, FCLM, FRLM, and CFRLM models.}
    \label{fig:preliminary_models}
\end{figure}

We depict the characteristics and differences of node-based (including the SCLM, MCLM, and p-median) and flow-based (including the FCLM, FRLM, and CFRLM) EVCI planning strategies in Fig. \ref{fig:preliminary_models}.

Besides the aforementioned optimization formulation for EVCI deployment, fuzzy decision-making, considering environmental, economic, and social factors, was introduced in \cite{guo2015_fuzzyTOPSIS,liu2019_multiCriteriaDecisionMaking,deb2022_MCLM_CSOandTLBO} for the charging station siting problem, providing well-rounded policy suggestions for city planners from more macro perspectives.

\subsubsection{Modeling Uncertain Routing Choices of EV Traffic Flows} \label{subsubsec:methodology_preliminaryModel_uncertainEVFlowModeling}

Previous studies have often overlooked the uncertainty of EV traffic flows, relying on the assumption that all EV drivers choose the shortest path to complete their O-D trips, which are calculated via the Dijkstra and Floyd search algorithms. However, this assumption may significantly deviate from realistic routing choices of EV drivers. Each O-D trip presents multiple viable routes, influenced by uncertain and uncontrollable individual driving behaviors, as well as various factors, such as traffic congestion and the limitation of driving range. To address this issue, He \textit{et al.}~\cite{he2015_MCLM_GA} introduced the deterministic user equilibrium (DUE) traffic assignment model to monitor actual EV routing decisions. The DUE model adheres to the Wardrop's equilibrium law~\cite{wardrop1952_UE} that all travelers select their optimal routes in a manner that no traveler can reduce their travel time by unilaterally changing routes. The equilibrium is reached when all allocated routes have the same travel time, not exceeding the travel time on any unused route. Their work proposed a bi-level optimization framework for EVCI planning. The DUE model, incorporating driving range constraints, served as the lower-level problem, providing routing choices for the upper-level MCLM model. The EV travel time was calculated using the Bureau of Public Roads (BPR) function~\cite{BPRFunction1984} taking traffic congestion into consideration. The bi-level optimization problem was transformed into a single-level mixed integer programming problem via the Karush–Kuhn–Tucker conditions and was subsequently solved by the GA algorithm. The bi-level framework, enhanced by the DUE model, was further extended to the flow-based EVCI planning model, such as the FCLM by He \textit{et al.}~\cite{he2018_FCLM_branchAndBound}, and was validated on synthetic traffic networks, including the Nguyen-Dupuis and Sioux Falls networks. Similar efforts have also been made in recent studies \cite{ma2021_FCLM_branchAndBound,duan2022_onlySizing,wang2019_CFRLM_branchAndBound,ferro2022_FCLM}.

The DUE traffic assignment model, while useful, falls short in capturing the variability of individual driving behaviors, since it assumes that each EV driver has full knowledge of traffic conditions and accurate travel time functions on each traffic link of O-D pairs. This assumption indicates that EV drivers always choose the optimal route to a charging station with minimum travel time. However, in reality, EV drivers must make routing decisions based on their perception of travel time, often leading to deviations from the actual optimal routes. To account for the discrepancy and model the stochastic routing choices arising from imperfect traffic information, the stochastic user equilibrium (SUE) was proposed by Fisk~\cite{fisk1980} with corresponding adaptions of the Wardrop's law -- assuming that all travelers choose their optimal route, such that no traveler can improve their \textit{perceived} travel time by unilaterally changing routes. In the SUE model, the Logit model is adopted to calculate the probability of a particular route being perceived as optimal with the minimum travel time. The randomness of routing decisions is modeled to follow a Gumbel distribution with zero mean and unit standard deviation, describing the uncertainty of individual routing behaviors. In the literature, Riemann \textit{et al.}~\cite{riemann2015_FCLM_branchAndBound} and Wang \textit{et al.}~\cite{wang2018_FCLM_MOEAD} have incorporated the SUE into their FCLM models. Specifically, Wang \textit{et al.}~\cite{wang2018_FCLM_MOEAD} further enhanced their model by modifying the link travel time function, replacing the BPR function with a more efficient volume-delay function (referred to as the conical congestion function) to improve the accuracy of the captured traffic flow.

\subsubsection{EVCI Resilience Under Contingency Events}

Besides the previously discussed uncertainties in traffic routing, other factors such as climate changes and contingency events also present additional challenges, especially in terms of infrastructure resilience in designing efficient planning strategies for EVCI. However, these factors have received comparatively less attention yet are essential for robust and reliable EVCI deployment. Zhang \textit{et al.}~\cite{zhang2023climate} incorporated flood resilience into the strategic placement of EV charging stations. The developed framework seeks to optimize the placement of charging stations by considering three main objectives: maximizing charging convenience, minimizing the impact of flood hazards, and reducing the interference with existing charging infrastructure. The effectiveness of this framework is demonstrated through a case study conducted in the Waikiki area. Similarly, Purba \textit{et al.}~\cite{Darmawi2024emergency} and Li \textit{et al.}~\cite{li2022emergency} have explored evacuation plans for EVs in emergency scenarios, such as hurricanes.

\subsection{EVCI Planning Strategies in Coupled Traffic Network and Electric Grid} \label{subsec:methodology_strategyInCoupledNetworks}

In addition to the primary objectives of minimizing investment costs or maximizing the captured traffic flows in a given traffic network, EVCI planning must also consider the critical interconnection between the traffic network and the electric grid to ensure grid stability and minimize energy loss in the power distribution network. Various studies have addressed these grid-related aspects, including:
\begin{itemize}
    \item Sadhukhan \textit{et al.}~\cite{sadhukhan2022_onlyGrid} and Liu \textit{et al.}~\cite{liu2013_onlyGrid} formulated the grid connection of EVCI as a network loss model, seeking to minimize transmission losses of the distribution network, which was tested on the IEEE $33$-bus~\cite{sadhukhan2022_onlyGrid} and $123$-bus~\cite{liu2013_onlyGrid} radial distribution systems, respectively.

    \item Deb \textit{et al.}~\cite{deb2019_MCLM_CSOandTLBO,deb2022_MCLM_CSOandTLBO} focused on the stability of node voltages in the electric grid, proposing a voltage sensitivity factor to describe potential impacts on voltage fluctuations.

    \item Almutairi \textit{et al.}~\cite{almutairi2022_chargingPortfolioGridImpact} explored optimal EVCI portfolios, including level 1, 2, and 3 chargers, aiming to smooth the total network load for reliable electricity supply. They introduced a peak power index to evaluate the performance of candidate EVCI portfolios.

    \item Mukherjee \textit{et al.}~\cite{mukherjee2023_onlyGrid} developed a cost-effective approach to locating and sizing chargers within a medium voltage distribution network. Their study considered a combination of single-port and multi-port chargers for EVCI deployment, with the latter capable of serving multiple EVs simultaneously, potentially leading to cost savings. Additionally, they factored the bidirectional power flows from EV charging, allowing for potential arbitrage opportunities via vehicle-to-grid (V2G) technology.
    \item As the adoption of EVs continues to rise, Tao \textit{et al.}~\cite{tao2023_onlyElectricGrid} provided a feasible expansion planning solution for both transmission and distribution networks to accommodate the increasing EV loads.
\end{itemize}

Besides planning EVCI on either the traffic network (discussed in Section \ref{subsec:methodology_preliminaryModel}) or the electric grid as presented above, numerous studies have focused on the context of the coupled systems, aiming to balance the interests of all involved stakeholders, including EV drivers (for charging accessibility), the city planner (for minimizing infrastructure costs), EVCI operator (for stable revenue stream), and the electric grid operator (for safe grid operations). For instance, Wang \textit{et al.}~\cite{wang2013_FCLM_crossEntropy} utilized the FCLM model to maximize captured traffic flow considering power flow constraints in a $25$-node traffic network overlaid with an IEEE $33$-bus distribution network. Node-based models, e.g., the SCLM, were also employed by Sadeghibarzani \textit{et al.}~\cite{sadeghibarzani2014_SCLM_GA} and Zheng \textit{et al.}, where the former minimized the EVCI development cost and grid energy loss (with grid constraints approximated using historical data), while the latter maximized the net revenue of charging stations considering investment, operation, and maintenance costs, as well as power flow constraints in the IEEE $15$ and $43$-bus distribution networks. Yao \textit{et al.}~\cite{yao2014_FCLM_MOEAD}, Mao \textit{et al.}~\cite{mao2021_FCLM_crossEntropy}, and Xiao \textit{et al.}~\cite{xiao2023_FCLM_GA} presented multi-objective collaborative planning strategies based on the FCLM model, to maximize the captured traffic flow while simultaneously minimizing the overall cost of investment and energy losses. Similar multi-objective frameworks have been introduced in studies such as \cite{wang2018_FCLM_MOEAD,zhang2018_CFRLM_branchAndBound,wang2018_FCLM_robust,wu2023_FCLM_PSO,wang2019_CFRLM_branchAndBound,li2023_FRLM_branchAndBound,ferro2022_FCLM,ahmad2023_SCLM,li2023_MCLM_crossEntropy}, with corresponding improvements shown as follows.
\begin{itemize}
    
    \item Wang \textit{et al.}~\cite{wang2018_FCLM_MOEAD} introduced a multi-stage and multi-objective collaborative planning model in the coupled traffic and power distribution networks. Specifically, two objectives, consisting of 1) minimizing investment and operations costs of the electric grid and 2) maximizing the annually captured EV flow, are simultaneously optimized.

    \item Zhang \textit{et al.}~\cite{zhang2018_CFRLM_branchAndBound} addressed the siting and sizing problems of fast charging stations within interconnected networks. They considered a cost-effective planning model and employed the CFRLM to ensure capturing maximum traffic flows.

    \item Wang \textit{et al.}~\cite{wang2018_FCLM_robust} and Wu \textit{et al.}~\cite{wu2023_FCLM_PSO} both explored uncertainties in future trends of EV adoption, comprising the growth rates of EV uptake and the distributions of EV charging loads in the electric grid. These uncertainties are addressed by designing a robust EVCI planning strategy with objectives of investment cost minimization and charging demand satisfaction. In addition, the work by \cite{wu2023_FCLM_PSO} also considered carbon emissions during EVCI planning.

    \item Wang \textit{et al.}~\cite{wang2019_CFRLM_branchAndBound} and Li \textit{et al}~\cite{li2023_FRLM_branchAndBound} discussed potential expansion on existing traffic and distribution networks. In \cite{wang2019_CFRLM_branchAndBound}, an optimal expansion model is determined by taking the siting and sizing of new charging stations, newly-built traffic roads, and new transmission lines in the electric grid into account. Additionally, the work by \cite{li2023_FRLM_branchAndBound} also analyzed EV drivers' traveling and charging behaviors for a more well-rounded expansion strategy.

    \item Ferro \textit{et al.}~\cite{ferro2022_FCLM} proposed a bi-level optimization framework for EVCI planning. The lower level formulates the uncertainties of EV drivers' routing choices and charging demands. Such uncertainties are embedded into the upper-level problem, whose objective is to minimize the network costs in both traffic and distribution networks.

    \item Ahamd \textit{et al.}~\cite{ahmad2023_SCLM} optimized the interests of EVCI investors, EV drivers, and distribution network operators, while exploring the possibility of introducing renewable energy resources at charging stations to reduce energy stress in the electric grid.

    \item Li \textit{et al.}~\cite{li2023_MCLM_crossEntropy} proposed a comprehensive multi-objective framework for EVCI planning in a realistic 183-node traffic network overlaid with a 33-region distribution network. Four objectives were developed to find the optimal planning solution, including maximizing captured traffic flow, minimizing charging time costs, minimizing EV traveling distance, and enhancing grid reliability.
    
\end{itemize}

\begin{figure}[!t]
    \centering
    \includegraphics[width=\linewidth]{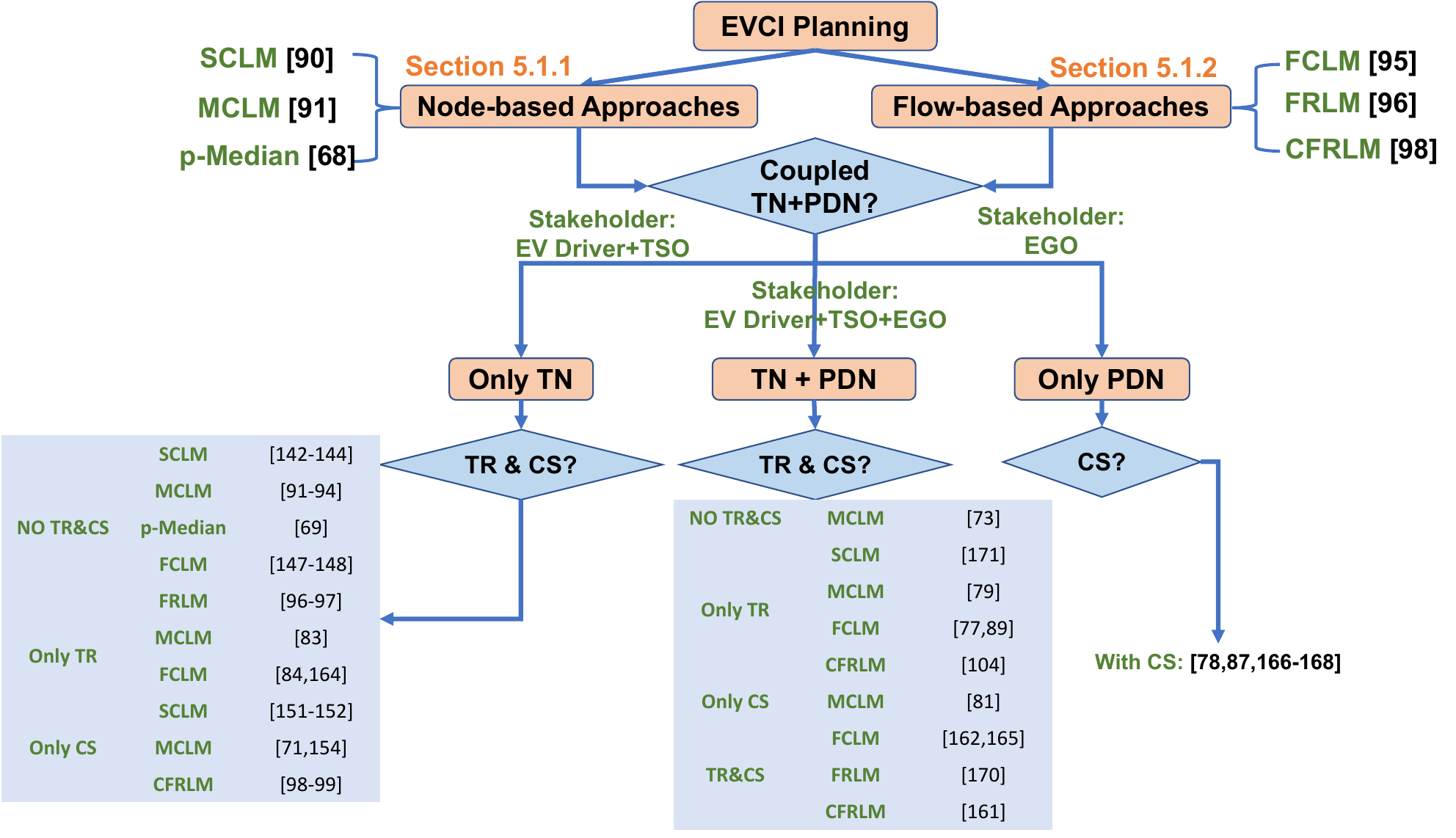}
    \caption{The hierarchical categories of reviewed solution methodologies. Here are explanations of abbreviations. TR: traffic routing; CS: capacity sizing; EGO: electric grid operator; TSO: transport system operator; TN: traffic network; PDN: power distribution network; SCLM: set covering location model; MCLM: maximum covering location model; FCLM: flow-capturing location model; FRLM: flow-refueling location model; CFRLM: capacitated flow-refueling location model.}
    \label{fig:hierarchical_methods}
\end{figure}

We summarize all reviewed studies, including key characteristics such as the employed facility location models in the traffic network (TN), or the power distribution network (PDN), or both, capacity sizing, routing choice modeling, objectives of EVCI planning, and algorithms to derive the optimal planning solutions, in Table \ref{tab:summary} for reference. Six objectives, including involved stakeholders, are presented with corresponding labels (from \textcircled{1} to \textcircled{6}) in Table \ref{tab:summary} and also shown as follows. -- \textcircled{1}: minimize infrastructure costs; \textcircled{2}: maximize covered EV charging demand; \textcircled{3}: minimize the EV weighted travel distance; \textcircled{4}: maximize captured EV traffic flows; \textcircled{5}: minimize the potential impacts of EV charging on the PDN; \textcircled{6}: maximize the operation revenues of EVCIs. These studies are also summarized in a hierarchical manner in Fig. \ref{fig:hierarchical_methods}. 

For simplicity, the adopted algorithms are shown with their abbreviations in Table \ref{tab:summary}. Their full forms are presented here: branch-and-bound (B\&B), particle swarm optimization (PSO), linear programming relaxation-based algorithm (LP-R), shared nearest neighbor clustering algorithm (SNN), genetic algorithm (GA), hybrid grey wolf optimization and particle swarm optimization (GWO-PSO), hybrid chicken swarm optimization and teaching-learning-based optimization (CSO-TLBO), multi-objective learner performance-based behavior algorithm (MOLPB), cross-entropy algorithm (CE), column and constraint algorithm (C\&C), multi-objective evolutionary algorithm based on decomposition (MOEA/D), data envelopment analysis (DEA), second-order cone programming (SOCP), non-dominated sorting genetic algorithm-II (NSGA-II), and modified primal-dual interior point algorithm (MPDIPA). 

In the above algorithms, B\&B, C\&C, and the interior point algorithm have gained widespread adoption in solving mixed integer linear programming problems. B\&B, in particular, serves as the default algorithm in several commercial optimization solvers such as Gurobi. However, with the expanding scale of both traffic networks and electric grids, leveraging these algorithms to derive an optimal planning solution is computationally expensive. Moreover, as planning strategies often consider multiple stakeholders with diverse corresponding objectives, these algorithms necessitate transforming the multi-objective problem into a weighted-sum single objective. Determining the appropriate weights for different objectives poses practical challenges. Therefore, heuristic algorithms, such as PSO, GA, and NSGA-II, stand out for their significantly lower computation burdens in searching the Pareto frontiers, i.e., the planning solutions under multi-objective frameworks.

\begin{center}
\begin{longtable}{c|c|c|c|c|c}
    \caption{Summary of reviewed studies.} \label{tab:summary} \\ 
         Ref. & \makecell{Location Model \\ in TN or PDN} & \makecell{Capacity \\ Sizing} & \makecell{Traffic \\ Routing} & \makecell{Objectives} & \makecell{Algorithm} \\

        \hline 
 
         \cite{toregas1971_SCLM} & SCLM+TN & \XSolidBrush & \XSolidBrush & \textcircled{1} & B\&B \\
         \cite{wang2009_SCLM_branchAndBound} & SCLM+TN & \XSolidBrush & \XSolidBrush & \textcircled{1} & B\&B \\
         \cite{wang2010_SCLM_branchAndBound} & SCLM+TN & \XSolidBrush & \XSolidBrush & \textcircled{1} & B\&B \\
         \cite{liu2012_SCLM_APSO} & SCLM+TN & \XSolidBrush & \XSolidBrush & \textcircled{1} & PSO \\
         \cite{jia2012_SCLM_branchAndBound} & SCLM+TN & \XSolidBrush & \XSolidBrush & \textcircled{1} & B\&B \\
         \cite{vazifeh2019_SCLM_GA} & SCLM+TN & \XSolidBrush & \XSolidBrush & \textcircled{1} & GA \\
         \cite{wang2017_SCLM} & SCLM+TN & \XSolidBrush & \XSolidBrush & \textcircled{1} & LP-R \\
         \cite{davidov2017_SCLM_branchAndBound} & SCLM+TN & \XSolidBrush & \XSolidBrush & \textcircled{1} & B\&B \\
         \cite{dong2016_SCLM_SNN} & SCLM+TN & \CheckmarkBold & \XSolidBrush & \textcircled{1} & SNN \\
         \cite{rajabi2017_SCLM_GA} & SCLM+TN & \CheckmarkBold & \XSolidBrush & \textcircled{1} & GA \\
         \cite{zhang2016_SCLM_PSO} & SCLM+TN & \CheckmarkBold & \XSolidBrush & \textcircled{1} & PSO \\
         \cite{sadeghibarzani2014_SCLM_GA} & SCLM+TN & \CheckmarkBold & \XSolidBrush & \textcircled{1}\textcircled{5} & GA \\
         \cite{ahmad2023_SCLM} & SCLM+TN\&PDN & \CheckmarkBold & \XSolidBrush & \textcircled{1}\textcircled{1}\textcircled{5}\textcircled{6} & GWO-PSO\\

         \hline\hline
         \cite{church1974_MCLM} & MCLM+TN & \XSolidBrush & \XSolidBrush & \textcircled{2} & B\&B \\
         \cite{xi2013_MCLM_branchAndBound} & MCLM+TN & \XSolidBrush & \XSolidBrush & \textcircled{2} & B\&B \\
         \cite{cai2014_MCLM} & MCLM+TN & \XSolidBrush & \XSolidBrush & \textcircled{2} & B\&B \\
         \cite{dong2014_MCLM_GA} & MCLM+TN & \XSolidBrush & \XSolidBrush & \textcircled{2} & GA \\
         \cite{jung2014_MCLM_GA} & MCLM+TN & \XSolidBrush & \XSolidBrush & \textcircled{2} & GA \\
         \cite{shahraki2015_MCLM_branchAndBound} & MCLM+TN & \XSolidBrush & \XSolidBrush & \textcircled{2} & B\&B \\
         \cite{tu2016_MCLM_GA} & MCLM+TN & \XSolidBrush & \XSolidBrush & \textcircled{2} & GA \\
         \cite{asamer2016_MCLM_branchAndBound} & MCLM+TN & \XSolidBrush & \XSolidBrush & \textcircled{2} & B\&B \\
         \cite{he2015_MCLM_GA} & MCLM+TN & \XSolidBrush & \CheckmarkBold & \textcircled{2} & GA \\
         \cite{cavadas2015_MCLM_pMedian_branchAndBound} & MCLM+TN & \CheckmarkBold & \XSolidBrush & \textcircled{2} & B\&B \\
         \cite{yang2017_MCLM_branchAndBound} & MCLM+TN & \CheckmarkBold & \XSolidBrush & \textcircled{2} & B\&B \\
         \cite{deb2022_MCLM_CSOandTLBO} & MCLM+TN\&PDN & \XSolidBrush & \CheckmarkBold & \textcircled{2}\textcircled{5} & CSO-TLBO \\
         \cite{pal2023_MCLM} & MCLM+TN\&PDN & \CheckmarkBold & \XSolidBrush & \textcircled{2}\textcircled{5} & MOLPB \\
         \cite{li2023_MCLM_crossEntropy} & MCLM+TN\&PDN & \XSolidBrush & \XSolidBrush & \textcircled{2}\textcircled{3}\textcircled{5} & CE \\

         \hline\hline 
         \cite{xu2013_pMedian_BPSO} & p-median+TN & \XSolidBrush & \XSolidBrush & \textcircled{3} & PSO \\
         \cite{an2014_pMedian_columnAndConstraint} & p-median+TN & \XSolidBrush & \XSolidBrush & \textcircled{3} &  C\&C \\
         \cite{he2016_pMedian_branchAndBound} & p-median+TN & \XSolidBrush & \XSolidBrush & \textcircled{3} & B\&B \\

         \hline\hline 
         \cite{hodgson1990_FCLM} & FCLM+TN & \XSolidBrush & \XSolidBrush & \textcircled{4} & B\&B \\
         \cite{wang2013_FCLM_branchAndBound} & FCLM+TN & \XSolidBrush & \XSolidBrush & \textcircled{4} & B\&B \\
         \cite{huang2015_FCLM_branchAndBound} & FCLM+TN & \XSolidBrush & \XSolidBrush & \textcircled{4} & B\&B \\
         \cite{sun2020_FCLMfast_MCLMslow_branchAndBound} & FCLM+TN & \XSolidBrush & \XSolidBrush & \textcircled{4} & B\&B \\
         \cite{ma2021_FCLM_branchAndBound} & FCLM+TN & \XSolidBrush & \CheckmarkBold & \textcircled{3}\textcircled{4} & B\&B \\
         \cite{riemann2015_FCLM_branchAndBound} & FCLM+TN & \XSolidBrush & \CheckmarkBold & \textcircled{4} & B\&B \\
         \cite{xiang2016_FCLM_branchAndBound} & FCLM+TN\&PDN & \CheckmarkBold & \XSolidBrush & \textcircled{1}\textcircled{4}\textcircled{5} & B\&B \\
         \cite{ferro2022_FCLM} & FCLM+TN\&PDN & \CheckmarkBold & \CheckmarkBold & \textcircled{1}\textcircled{4}\textcircled{5}\textcircled{6} & B\&B \\
         \cite{wang2018_FCLM_MOEAD} & FCLM+TN\&PDN & \CheckmarkBold & \CheckmarkBold & \textcircled{1}\textcircled{4}\textcircled{5}\textcircled{6} & MOEA/D \\
         \cite{mao2021_FCLM_crossEntropy} & FCLM+TN\&PDN & \CheckmarkBold & \XSolidBrush & \textcircled{1}\textcircled{4}\textcircled{5} & CE \\
         \cite{xiao2023_FCLM_GA} & FCLM+TN\&PDN & \CheckmarkBold & \XSolidBrush & \textcircled{1}\textcircled{4}\textcircled{5} & GA \\
         \cite{wang2018_FCLM_robust} & FCLM+TN\&PDN & \CheckmarkBold & \XSolidBrush & \textcircled{1}\textcircled{4}\textcircled{5} & DEA \\
         \cite{wu2023_FCLM_PSO} & FCLM+TN\&PDN & \CheckmarkBold & \XSolidBrush & \textcircled{1}\textcircled{4}\textcircled{5} & PSO \\
         
         \hline\hline
         \cite{kuby2005_FRLM} & FRLM+TN & \XSolidBrush & \XSolidBrush & \textcircled{4} & B\&B \\
         \cite{kim2012_FRLM_branchAndBound} & FRLM+TN & \XSolidBrush & \XSolidBrush & \textcircled{4} & B\&B \\
         \cite{wu2017_FRLM_branchAndBound} & FRLM+TN & \XSolidBrush & \XSolidBrush & \textcircled{4} & B\&B \\
         \cite{li2023_FRLM_branchAndBound} & FRLM+TN\&PDN & \CheckmarkBold & \CheckmarkBold & \textcircled{1}\textcircled{4}\textcircled{5} & B\&B \\

         \hline\hline
         \cite{upchurch2009_CFRLM} & CFRLM+TN & \CheckmarkBold & \XSolidBrush & \textcircled{4} & B\&B \\
         \cite{li2016_CFRLM_GA} & CFRLM+TN & \CheckmarkBold & \XSolidBrush & \textcircled{4} & GA \\
         \cite{zhang2018_CFRLM_secondOrderCone} & CFRLM+TN & \CheckmarkBold & \XSolidBrush & \textcircled{4} & SOCP \\
         \cite{wang2018_CFRLM_GA} & CFRLM+TN & \CheckmarkBold & \XSolidBrush & \textcircled{4} & GA \\
         \cite{wang2019_CFRLM_branchAndBound} & CFRLM+TN\&PDN & \CheckmarkBold & \CheckmarkBold & \textcircled{1}\textcircled{4}\textcircled{5} & B\&B \\
         \cite{zhang2018_CFRLM_branchAndBound} & CFRLM+TN\&PDN & \CheckmarkBold & \XSolidBrush & \textcircled{1}\textcircled{4}\textcircled{5} & B\&B \\

         \hline\hline
         \cite{sadhukhan2022_onlyGrid} & Only PDN & \CheckmarkBold & \XSolidBrush & \textcircled{5} &  NSGA-II \\
         \cite{liu2013_onlyGrid} & Only PDN & \CheckmarkBold & \XSolidBrush & \textcircled{5} &  MPDIPA \\
         \cite{almutairi2022_chargingPortfolioGridImpact} & Only PDN & \CheckmarkBold & \XSolidBrush & \textcircled{5} &  B\&B \\
         \cite{mukherjee2023_onlyGrid} & Only PDN & \CheckmarkBold & \XSolidBrush & \textcircled{5}\textcircled{6} & B\&B \\
         \cite{tao2023_onlyElectricGrid} & Only PDN & \CheckmarkBold & \XSolidBrush & \textcircled{5}\textcircled{6} & B\&B \\

\end{longtable}
\end{center}

Apart from the above literature focusing on strategic EVCI spatial placement, designing appropriate investment plans, as well as effective budget allocation, during the planning horizon are also essenital for an efficient EVCI planning strategy. Rehman \textit{et al.}~\cite{Rehman2024} explored strategic multi-period coordinated planning for optimally siting and sizing fast charging stations in a highway transportation and power distribution network. In particular, this multi-period planning approach recognizes the dynamic nature of technology costs, deciding the number and location of charging stations over three distinct periods to cover the entire planning horizon. The model outlines the importance of not only investing immediately but also planning for future expansions and technological upgrades. This is contrasted with the traditional forward-myopic method, which handles planning in a sequential single-period manner, potentially missing efficiencies from a more integrated long-term perspective. Their model incorporates spatio-temporal variations in EV charging demand, which influence both the timing and scaling of infrastructure investment. Similarly, Borozan \textit{et al.}~\cite{borozan2022} integrated EV smart charging into network expansion planning, highlighting the need for investments in EVCI to be timed with considerations of long-term uncertainties in electric grids. Besides, a whitepaper from the European Automobile Manufacturers Association (ACEA)~\cite{eu2022_masterPlan} stresses that investments should be strategically allocated to create a balanced and accessible charging network across all regions of Europe. This includes prioritizing both urban and rural areas to ensure comprehensive coverage.

\section{Research gaps and call for contributions} \label{sec:reserach_gaps}

\begin{figure}[!t]
    \centering
    \includegraphics[width=\linewidth]{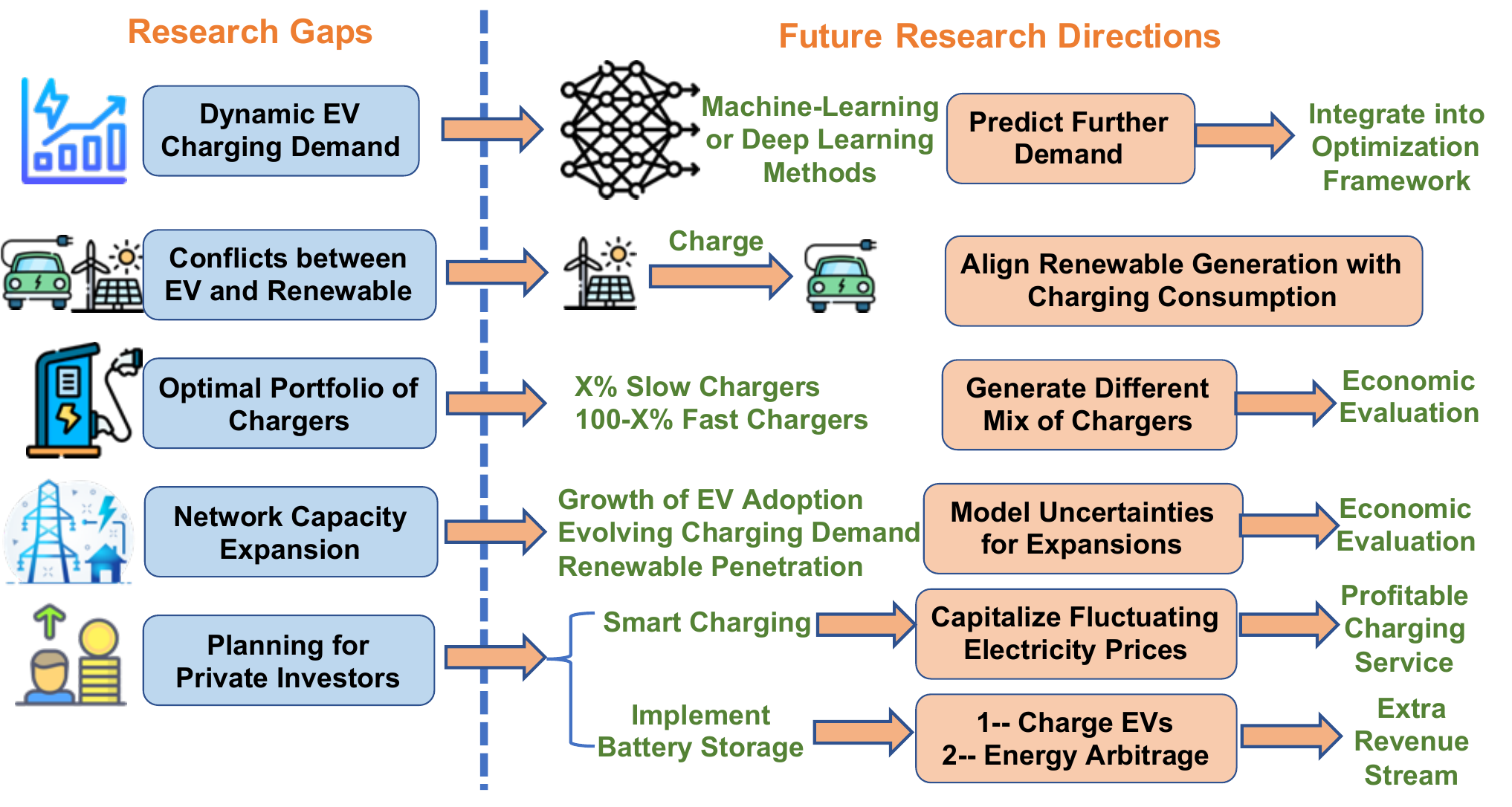}
    \caption{The illustration of research gaps and corresponding research directions.}
    \label{fig:research_gap_illustration}
\end{figure}

As discussed above, previous studies have significantly progressed in designing efficient EVCI planning strategies. However, several research gaps still exist, requiring further contributions to this field. The identified gaps are pointed out as follows, each of which is followed by possible solution methods to address them for more efficient planning strategies. We also depict both research gaps and solutions in Fig. \ref{fig:research_gap_illustration}.

\begin{itemize}
    \item \textit{Incorporating Dynamic EV Charging Demand}: While the stochastic routing choices can be modeled using traffic assignment models such as the DUE and SUE, current literature often assumes stationary EV charging demands based on historical data. With the expected rapid growth in EV adoption, integrating dynamic and time-varying EV charging demands into the EVCI planning problem is essential to create efficient and cost-effective solutions that can cater to the evolving needs of EV drivers without necessitating additional infrastructure expansion. 

    \textbf{Possible Solution:} Rather than relying solely on a static analysis of EV charging demand derived from historical data, harnessing the power of emerging machine-learning or deep-learning techniques to predict future charging demand can significantly enhance the effectiveness of planning solutions, especially in the context of the rapidly accelerating EV transition. These predictions can seamlessly integrate into the optimization framework, thereby facilitating the generation of more practical and timely EVCI planning solutions.

    \item \textit{Coordinating Vehicle Electrification with Grid Decarbonization}: Electrifying road transport is a vital step towards decarbonization. However, the electrification of vehicles can challenge the net-zero emission goals of the electric grid. The main reason lies in the underutilization of renewable energy resources (such as wind and solar energy). Although such distributed energy resources have been expanded into the distribution network, the current EV charging behaviors still require coal or gas generators to ensure sufficient generation capacity supporting the stochastic and surging load of EV fleets or large-scale residential (or workplace) charging. Consequently, the uncertain EV charging demand may create barriers to further renewable penetration in the grid. Hence, an efficient EVCI planning strategy should consider utilizing abundant renewables to charge EVs rather than requiring more thermal generator participation, aligning with the decarbonization and net-zero targets in the coupled systems. Moreover, the impacts of both increasing renewables and EV charging demand on the grid should also be examined to ensure safe grid operations and reliable electricity supply.

    \textbf{Possible Solution:} To effectively coordinate EV adoption within the traffic network and the increasing renewable penetration in the electric grid, a mutually beneficial scenario can be achieved by aligning renewable generation with EV charging consumption. This synergy is best realized through the implementation of smart charging operations, which involve dynamically allocating charging demand and optimizing the utilization of electricity supply from renewable sources. By doing so, renewable energy can be efficiently consumed by EV charging, mitigating the need for curtailment during periods of oversupply. This approach not only encourages the widespread adoption of renewable energy resources but also establishes a win-win situation for all stakeholders involved. As the grid undergoes a transition with the gradual exit of thermal generators and the expansion of renewable generators, the increased integration of renewables into the grid contributes to a reduction in electricity prices. Consequently, this development enhances the cost-effectiveness of EV charging for drivers, providing them with an additional incentive to adopt electric vehicles. In essence, the collaborative efforts between renewable energy integration and smart charging operations contribute to a sustainable and economically viable evolution in the transportation and energy sectors.

    \item \textit{Optimal Portfolio of Charging Technologies}: Evaluating the optimal portfolio of various charging technologies (i.e., level 1, 2, and 3 chargers) in diverse geospatial locations (e.g., residential areas, workplaces, or highways) is 
    pivotal to meeting individual charging needs without excessive investment. However, most literature only focused on a single type of charging station (e.g., fast charging). Though few considered placing various types of chargers throughout the traffic network, they still assumed that one implemented charging station is only equipped with one specific kind of charger. Therefore, there is a strong need for comprehensive studies considering diverse charger types and locations.

    \textbf{Possible Solution:} In the process of determining the capacity for each charging station, it is imperative to integrate a well-defined optimization framework that parameterizes the composition of chargers. This involves specifying the percentages of different charger types within a single charging station, a consideration particularly crucial in diverse contexts such as commercial or residential areas. Additionally, tailoring the capacity sizing strategy to the unique characteristics of different areas is essential, given the substantial variations in charging behaviors and requirements. To achieve an economically optimized mix of charging technologies, a thorough examination of charging dynamics in various settings is warranted. This economic analysis should be instrumental in identifying the most suitable and efficient blend of charging infrastructure for each location. In essence, a context-specific and economically grounded approach to capacity sizing ensures the optimal deployment of charging stations, aligning them with the distinct needs and patterns of charging demand in different areas.

    \item \textit{Capacity Expansion of Existing Charging Resources}: Besides building new EVCI (which most of the research focused on), a cost-effective EVCI planning strategy should account for the potential capacity expansion of existing EV charging resources, alongside the construction of EVCI in new candidate sites in traffic networks. Such a planning strategy optimizes the utilization of current assets, thereby reducing the financial burden associated with new EVCI construction.

    \textbf{Possible Solution:} The development of expansion strategies necessitates a proactive approach that involves modeling uncertainties surrounding future information. Key variables include the growth rates of EV adoption, the dynamic evolution of EV charging demand, and the anticipated integration of renewable energy into the electric grid. To effectively address these uncertainties, it is essential to employ scenario modeling, examining diverse potential outcomes. Within the framework of these modeled scenarios, the determination of optimal expansion strategies hinges on economic evaluations. By assessing the economic viability of different expansion paths against the backdrop of current network structures, decision-makers can strategically navigate the complexities associated with the evolving landscape of EV adoption and renewable energy integration. This forward-looking and economically grounded approach ensures that expansion efforts are not only resilient to uncertainties but also aligned with the most favorable economic outcomes under various envisioned scenarios.

    \item \textit{EVCI Planning Strategies for Private Charging Station Operators}: The perspective of third-party charging station operators, often private market participants, remains underrepresented in EVCI planning studies. These operators aim to maximize profits and compete with other service providers. Consequently, there is also a need for research that addresses siting, sizing, and capacity expansion from the viewpoint of these market players, not solely from the city planner's perspective.

    \textbf{Possible Solution:} Private investors can capitalize on the dynamic nature of electricity prices, leveraging fluctuations to enhance the profitability of charging services. With the rising adoption of renewable energy furthering the volatility of the electricity market, an attractive strategy involves integrating a battery energy storage system (BESS) alongside charging stations. This co-location offers a dual benefit -- 1) the BESS can be charged during periods of abundant renewable generation and lower daytime electricity prices. The stored energy can be released for charging customers in cases of higher electricity prices; 2) the BESS can also actively participate in the electricity market, performing arbitrage operations to create an additional revenue stream.

\end{itemize}

\section{Conclusion} \label{sec:conclusion}
In this review paper, we have examined the current state of EVCI planning. Our exploration begins with an overview of the contemporary efforts in EVCI planning across key countries in the global EV market. This analysis reveals a notable lag in the development of EVCI planning compared to the escalating adoption of EVs. Three significant barriers hindering the widespread implementation of EVCI are identified: inadequate EVCI charging services, low utilization rates of public EVCI, and the intricate integration of EVCI into the electric grid. 

To delve deeper into the intricacies of the infrastructure planning problem, we introduce stakeholders and various charging technologies applicable to EVCI. Stakeholders, including EV drivers, transport system operators, and electric grid operators, exhibit distinct interests in accessibility, infrastructure costs, and distribution network stability. Charging technologies, namely level 1, 2, and 3 chargers, are explored, with considerations for their potential applications in diverse geographical areas based on varying charging speed requirements. We further categorize EV charging into residential, destination, and en-route charging, aligning with level 1, a combination of level 2 and 3, and level 3 charging in practical scenarios. The paper proceeds to a thorough review of EVCI planning solution methodologies, primarily encompassing node-based and flow-based approaches. Node-based approaches concentrate on specific traffic nodes for charging demand, while flow-based approaches articulate traffic flows through O-D trips. Addressing key challenges, such as capacity sizing of charging stations, uncertainty in EV drivers' routing choices, and integration of EVCI into the electric grid, we present solutions proposed in the existing literature. 

In conclusion, we underscore the need for further contributions to bridge existing research gaps. Key areas for advancement include: 1) incorporating evolving and time-varying charging demand into EVCI planning solutions, 2) coordinating vehicle electrification with grid decarbonization, 3) exploring the optimal portfolio of charging technologies, 4) considering additional capacity expansion of charging resources, and 5) examining the economic potential of operating charging stations from the perspective of private EVCI operators. This call for continued research aims to propel the evolution of EVCI planning and ensure its alignment with the dynamic landscape of electric mobility.

\section*{Data availability statement}
No new data were created or analyzed in this study.

\section*{Acknowledgements}
This work is supported in part by the Australian Research Council (ARC) Discovery Early Career Researcher Award (DECRA) under Grant DE230100046.

\section*{References}

\bibliographystyle{ieeetr}
\bibliography{ref}

\end{document}